\documentclass[12pt,oneside,english]{amsart}
\usepackage[T1]{fontenc}
\usepackage[latin1]{inputenc}
\usepackage{geometry}
\usepackage{amssymb}

\makeatletter
\theoremstyle{plain}
\newtheorem{thm}{Theorem}[section]
\numberwithin{equation}{section} %% Comment out for sequentially-numbered
\numberwithin{figure}{section} %% Comment out for sequentially-numbered
\theoremstyle{plain}
\theoremstyle{remark}
\newtheorem*{acknowledgement*}{Acknowledgement}
\theoremstyle{plain}
\newtheorem{prop}[thm]{Proposition} %%Delete [thm] to re-start numbering
\theoremstyle{plain}
\newtheorem{lem}[thm]{Lemma} %%Delete [thm] to re-start numbering
\theoremstyle{definition}
\newtheorem{defn}[thm]{Definition}
\theoremstyle{remark}
\newtheorem{notation}[thm]{Notation}
\theoremstyle{remark}
\newtheorem{rem}[thm]{Remark}
\theoremstyle{plain}
\newtheorem{cor}[thm]{Corollary} %%Delete [thm] to re-start numbering
\theoremstyle{remark}

\theoremstyle{remark}

\newcommand{\poly}[2]{\ensuremath{\mathbb{#1}\langle{#2}\rangle}}
\newcommand{\polyxk}{\mathbb{R}\langle{x_1,\ldots,x_k}\rangle}
\DeclareMathOperator{\Tr}{Tr}

\usepackage{babel}
\makeatother
\begin{document}

\title{Information Geometry of Random Matrix Models}

\author{Dan Shiber}

\thanks{The author was supported in part by an NSF VIGRE fellowship.}

\begin{abstract}
In this paper we develop the theory of information geometry for
single random matrix models, with two goals: proving a Cramer-Rao
theorem for estimators on random matrices, and calculating the
Legendre transform of pressure and entropy with respect to a metric
duality. Consequently, in the large $n$ limit we recover several quantities from free
probability: Voiculescu's conjugate variable is the tangent vector
to the GUE perturbation model, giving rise to a metric which turns
out to be the free Fisher information measure; Hiai's Legendre transform of
free pressure agrees with our Legendre transform of pressure; and
Speicher's covariance of fluctuations naturally arises as the metric
on the random matrix model obtained from the fluctuation functions.
\end{abstract}

\maketitle

~

\subsection{Introduction\label{sec:Introduction}}

~

Inspired by the work of \cite{AmNa}, we treat random $n\times{n}$ matrix models of the form\\ $\exp\left(-n\Tr\left(p(A)+\psi(n)\right)\right)$ with $p\in\mathbb{R}\langle{x}\rangle$ and $\psi(n)=\frac{1}{n^2}\log\int\exp\left(-n\Tr\left(p(A)\right)\right)dA$ as statistical models and construct their information geometry. This achieves two goals: it proves the Cramer-Rao theorem, which is a Cauchy-Schwartz inequality on polynomial functions of the random matrix (Section \ref{sec:Cramer-Rao-Theorem}); and it calculates the entropy as the Legendre transform of pressure (Section \ref{sec:Legendre-Transform}).

In Section \ref{sec:The-n-goes-to-infinity-case} we relate our construction to free probability by considering the limit as the matrix size $n$ approaches infinity. We show that the information geometric quantities converge. The pressure, entropy, and Legendre transform converge to the free pressure, free entropy, and free Legendre transform \cite{Hi} respectively. We also show that the information geometry of a Gaussian perturbation model converges to the free Fisher information measure \cite{Voi1}. Finally we note the relation to the Free Cramer-Rao Theorem \cite{Voi2} and the fluctuations of random matrices \cite{MSp}.

Following is a quick review of classical information geometry, meant as a motivation for our development. The familiar reader may skip along to Section \ref{sec:Basic-Notions}.

~

\begin{acknowledgement*}
I would like to thank my advisor, Dimitri Shlyakhtenko, for
suggesting this project and the countless discussions which led to
its fruition.
\end{acknowledgement*}

~

\subsection{Review of Classical Information Geometry}

~

Classical information geometry may be viewed as the standard framework for
doing convex analysis (finding minima/maxima) on real-valued
functions of random variables. Given a random variable $X_{\theta}$
whose distribution function belongs to a parametric model
$\left\{\left.q_{\theta}(x)dx\right|\theta\in\Theta\subset\mathbb{R}^m\right\}$,
and functions (estimators) $\xi_1,\ldots,\xi_m\in C(\mathbb{R})$,
one is interested in measuring the sensitivity of
$\xi_{1}(X_{\theta}),\ldots,\xi_{m}(X_{\theta})$ to changes in
$\theta$. This analysis is done following the presentation of
\cite{AmNa} using the methods of differential geometry, and the
resulting theorem is a lower bound on the covariance of the
deviations of the $\xi_{i}$'s as follows.

~

A \emph{statistical model} $S$ is a family of probability
distributions on $\mathbb{R}$ parameterized by finitely many real
parameters, $S=\left\{ q_{\theta}(x)dx\left|\theta\in\Theta,\
q_{\theta}(x)>0,\ \int q_\theta(x)dx=1\right.\right\}$ with
$\Theta\subset\mathbb{R}^{m}$ open.

~

An \emph{exponential family} is a statistical model with
$q_{\theta}(x)=\exp\left(p(x)+\sum\theta_{i}f_{i}(x)+\psi(\theta)\right)$,
where
$\psi(\theta)=-\log\int\exp\left(p(x)+\sum\theta_{i}f_{i}(x)\right)dx$
and $p,f_1,\ldots,f_m\in C(\mathbb{R})$ such that $\psi(\theta)$
converges. We denote
$p_\theta(x)=p(x)+\sum\theta_{i}f_{i}(x)+\psi(\theta)$.

~

An exponential family $S$ is a manifold under the map
$\exp\left(p_\theta(x)\right)\mapsto\theta$. Its tangent space is
the vector space of random variables
\begin{align*}T_\theta{S}=\textrm{span}\left\{\left.\frac{\partial}{\partial\theta_{i}}\exp\left(p_\theta(x)\right)(X_{\theta})\right|X_{\theta}\sim
\exp\left(p_\theta(x)\right)\right\}_{i=1}^{m}\mathop{\simeq}^{\log}
\textrm{span}\left\{\left.\frac{\partial}{\partial\theta_{i}}p_{\theta}(X_{\theta})\right|X_{\theta}\sim
\exp\left(p_\theta(x)\right)\right\} _{i=1}^{m}.\end{align*}

~

There is a natural $L^{2}$-structure on this space which allows us
to define an inner-product
\begin{align*}\left\langle f,g\right\rangle_{\theta}=\int f(x)g(x)\exp\left(p_\theta(x)\right)dx,\ f,g\in{T}_\theta{S},\end{align*}
and this gives the Fisher Information Metric
\begin{align*}g_{ij}(\theta)=\left\langle\frac{\partial}{\partial\theta_{i}}p_{\theta},\frac{\partial}{\partial\theta_{j}}p_{\theta}\right\rangle_{\theta}.\end{align*}
The $L^{2}$-structure also identifies the \emph{potential}
$p_\theta$ as $-\int_{-\infty}^{x}d_{\theta}^{*}(1)(y)dy$, where
\begin{align*}d_{\theta}:L^{2}(\exp\left(p_\theta(x)\right))\rightarrow
L^{2}(\exp\left(p_\theta(x)\right)),\ f\mapsto f^\prime\end{align*}
is the unbounded differentiation operator. This follows from the
calculation
\begin{align*}
&\big\langle{d_{\theta}^{*}(1)(x),h(x)}\big\rangle_{L^{2}\left(\exp(p_{\theta})\right)}=\int1\cdot{h}^{\prime}(x)\exp\left(p_{\theta}(x)\right)dx=\\ &-\int{h}(x)p_{\theta}^{\prime}(x)\exp\left(p_\theta(x)\right)dx=\left\langle{h(x),-p_{\theta}^{\prime}(x)}\right\rangle_{L^2\left(\exp\left(p_\theta(x)\right)dx\right)},
\end{align*}
for any polynomial $h$. Thus, the tangent space consists of partial
derivatives of $-\int_{-\infty}^{x}d_{\theta}^{*}(1)(y)dy$ with
respect to $\theta_i$.

~

In this framework \cite{AmNa} prove the Cramer-Rao theorem:

\begin{thm}
Let $S$ be an exponential family with Fisher information metric $g$,
and let \\
$\xi_{1},\ldots,\xi_{m}:\mathbb{R}\rightarrow\mathbb{R}^m$ be
\emph{unbiased estimators} i.e.
$\int\xi_{i}(x)\exp\left(p_\theta(x)\right)dx=\theta_i$.  Then
\begin{align*}
\left\langle\xi_{i}(X_{\theta})-\theta_i,\xi_{j}(X_{\theta})-\theta_j\right\rangle\geq{g}_{ij}^{-1}(\theta)
\end{align*}
in the sense of positive semi-definite matrices.
\end{thm}

~

Next, to find minima/maxima on $S$ and calculate the Legendre
transform of $\psi$, \cite{AmNa} specify a second derivative (the
tangent space being the first derivative).  This is done by fixing
an affine connection (see \cite{BiGo}), which is given in
coordinates by
\begin{align*}
\Gamma_{ijk}^{(\alpha)}(\theta)=\int\frac{\partial}{\partial\theta_{k}}p_{\theta}(x)\cdot\frac{\partial^{2}}{\partial\theta_{i}\partial\theta_{j}}p_{\theta}(x)\cdot{p}_{\theta}(x)dx+\frac{1-\alpha}{2}\int\frac{\partial}{\partial\theta_{k}}p_{\theta}(x)\cdot\frac{\partial}{\partial\theta_{i}}p_{\theta}(x)\cdot\frac{\partial}{\partial\theta_{j}}p_{\theta}(x)\cdot{p}_{\theta}(x)dx
\end{align*}
where $\alpha$ is a parameter for the
amount of curvature. For example, an exponential family is flat for
$\alpha=1$, and a mixture family ($g_{\theta}=\mu_{\left(X+\theta
Y\right)/\left(1+\theta\right)}$ with $X$ and $Y$ independent) is
flat for $\alpha=-1$.

~

These connections are distinguished by the duality of the
$(\alpha)-$ and $(-\alpha)-$connections with respect to $g$:
\begin{prop}
$\frac{\partial}{\partial\theta_{k}}g_{ij}(\theta)=\Gamma_{kij}^{(\alpha)}(\theta)+\Gamma_{kji}^{(-\alpha)}(\theta)$.
\end{prop}

This allows \cite{AmNa} to prove
\begin{thm}
Let $S$ be a manifold with metric $g$, a pair of dual affine
connections $\Gamma,\Gamma^{*}$, and a smooth function
$f:M\rightarrow\mathbb{R}$. If $\theta^{\prime} \in M$ satisfies
$\frac{\partial}{\partial\theta_{i}}f(X_{\theta^{\prime}})=0$ and
$\frac{\partial}{\partial\Gamma(\theta)}f(X_{\theta^{\prime}})\geq
g$ in the sense of positive semi-definite matrices, where
$\frac{\partial}{\partial\Gamma(\theta)}$ is the covariant
derivative, then $\exists\tilde{\Theta}$ a small neighborhood of
$\theta^\prime$ such that $f(X_{\theta^\prime})=\sup_{\theta\in\tilde{\Theta}}f(X_{\theta})$.
\end{thm}

~

~

\section{Basic Notions\label{sec:Basic-Notions}}

\subsection{Manifold\label{sec:The-Manifold}}

~

We start with a random $n\times{n}$ self-adjoint matrix $A$ with complex entries, distributed according to
\begin{equation}
\exp\left(-n\Tr\left(p(A)+\psi(n)\right)\right)dA \textrm{ on
}M_{n}^{SA}(\mathbb{C})\label{eq:exponential}
\end{equation}
where $p\in\poly{R}{x}$ is convex, and $\psi(n)=\frac{1}{n^{2}}\log\int_{M_{n}^{SA}(\mathbb{C})}\exp\left(-n\textrm{Tr}\left(p(A)\right)\right)dA$ is the normalization to a probability measure.  The quantity $\psi$ is
also known as the pressure \cite{Hi}. In this paper
$\Tr:M_{n}^{SA}(\mathbb{C})\rightarrow\mathbb{C}$ by
$A\mapsto\sum_{i=1}^{n}A_{ii}$, and $dA=\prod_{\substack{1\leq i\leq
j\leq n}}d\mathfrak{Re}(A_{ij})d\mathfrak{Im}(A_{ij})\prod_{1\leq
i\leq n}dA_{ii}$. We recall a useful fact which guarantees
convergence of this model \cite{Bi}:

~

\begin{lem}\label{lem:biane}
Given $p\in{C}^{2}(\mathbb{R})$ convex, $\exists{!}q:\mathbb{R}\rightarrow\mathbb{R}$ Borel, $q(x)\geq0$, $\int_\mathbb{R}q(x)dx=1$ such that $\forall f\in C(\mathbb{R})$,
\begin{align*}\frac{1}{n}\int\Tr\left(f(A)\right)\exp\left(-n\Tr\left(p(A)+\psi(n)\right)\right)dA\rightarrow\int
f(x)q(x)dx\end{align*} where $q$ is defined by the integral equation
\begin{align*}
2\textrm{pr.v.}\int\frac{q(y)}{y-x}dy = p^{\prime}(x)\textrm{,
}\forall x\in\textrm{supp}(q).\end{align*}
\end{lem}

~

\begin{defn}
An \emph{exponential family} is a family of distributions on
$M_{n}^{SA}(\mathbb{C})$ of the form

\begin{align*}
S=\left\{
\left.\exp\left(-n\textrm{Tr}\left(p(A)+\sum_{i=1}^{m}\theta_{i}F_{i}(A)+\psi(\theta,n)\right)\right)\right|\theta\in\Theta\right\}
\textrm{,}\end{align*} where
\begin{align*}\psi(\theta,n)=\frac{1}{n^{2}}\log\int\exp\left(-n\textrm{Tr}\left(p(A)+\sum_{i=1}^{m}\theta_{i}F_{i}(A)\right)\right)dA\end{align*}
is the normalization constant, $F_{1},\ldots,F_{m}\in C(\mathbb{R})$
are the perturbation functions, and $\Theta\subset\mathbb{R}^{m}$ is
an open set of parameters (chosen so that the integral in the
definition of $\psi$ converges).
\end{defn}

\begin{notation}
Write $p_{\theta,n}(A)=p(A)+\sum_{i=1}^{m}\theta_{i}F_{i}(A)+\psi(\theta,n)$,
$p_{\theta}(A)=p(A)+\sum_{i=1}^{m}\theta_{i}F_{i}(A)$, and
$d\mu_{\theta,n}(A)=\exp\left(-n\textrm{Tr}\left(p_{\theta,n}(A)\right)\right)dA$.
\end{notation}
\begin{rem}
$S$ is a manifold under the chart \begin{align*}d\mu_{\theta,n}(A)\mapsto\theta.\end{align*}
\end{rem}

~

\subsection{Tangent Space\label{sec:The-Tangent-Space}}

~

In analogy to the classical case, we define the tangent space by
identifying the potential of the limit distribution $q_{\theta}$.
Since $q_\theta$ is the distribution of a noncommutative random variable,
the free difference quotient plays the role of the derivative:
\begin{align*}
\partial_\theta:L^2(\mathbb{R},q_\theta)\rightarrow{L}^2(\mathbb{R},q_\theta)\otimes{L}^2(\mathbb{R},q_\theta),\ f(x)\mapsto\frac{f(x)-f(y)}{x-y}
\end{align*}
(this operator is explained in \cite{Voi1}). It is a densely-defined
derivation with
$\textrm{domain}\left(\partial_{\theta}\right)=\textrm{polynomials}$.
For $g\in\textrm{domain}\left(\partial_{\theta}\right)$ we have
\begin{align*}
&\left\langle\partial_{\theta}^{*}\left(1\otimes1\right)(x),g(x)\right\rangle_{L^2\left(\mathbb{R},q_{\theta}\right)}=\left\langle1\otimes1,\left(\partial_{\theta}g\right)\left(x,y\right)\right\rangle_{L^2\left(\mathbb{R},q_{\theta}\right){\otimes}L^2\left(\mathbb{R},q_{\theta}\right)}=\\
&\int\frac{g(x)-g(y)}{x-y}q_{\theta}(x)q_{\theta}(y)dxdy=-2\int{g}(x)\left(\int\frac{q_{\theta}(y)}{y-x}dy\right)q_{\theta}(x)dx=\\
&\left\langle{p_{\theta}^{\prime}(x),g(x)}\right\rangle_{L^2\left(\mathbb{R},q_{\theta}\right)}.
\end{align*}
Therefore, the limit potential is identified as $\partial_\theta^*(1\otimes1)=\int_{-\infty}^{x}p_{\theta}^{\prime}(y)dy=p_{\theta}(x)$.

~

\begin{defn}
The potential of $d\mu_{\theta,n}(A)$ is $p_{\theta,n}(A)$.
\end{defn}
\begin{prop}
The tangent space is given by $T_{\theta}S=\textrm{span}\left\{ \left.\frac{\partial}{\partial\theta_{i}}p_{\theta,n}(A)\right|A\sim d\mu_{\theta,n}(A)\right\} _{i=1}^{m}$,
regarded as a vector space of random variables.
\end{prop}
\begin{proof}
Fix $\theta\in\Theta$. Each $\frac{\partial}{\partial\theta_1}p_{\theta,n},\ldots,\frac{\partial}{\partial\theta_m}p_{\theta,n}$ defines a curve through $\theta$:
\begin{align*}
\gamma_i(t)=\exp\left(-nTr\left(p_{\theta,n}(A)+t\frac{\partial}{\partial\theta_{i}}p_{\theta,n}(A)+\psi(t,n)\right)\right).
\end{align*}

~

Conversely, suppose $h\in\poly{R}{x}$ is convex and
$\exp\left(-n\Tr\left(h(A)+\phi(n)\right)\right)\in S$, with
$\phi(n)=\frac{1}{n^2}\log\int\exp(-n\Tr(h(A)))dA$.  Get
$\theta^{\prime}$ such that $\exp(-n\Tr(h(A)+\phi(n)))dA=d\mu_{\theta^{\prime},n}(A)$. Then by Lemma \ref{lem:biane}, $\exp(-n\Tr(h(A)+\phi(n))$ and
$d\mu_{\theta^{\prime},n}$ both converge to $q_{\theta^{\prime}}$
satisfying
\begin{align*}
h^\prime(x)=2\textrm{pr.v.}\int\frac{q_{\theta^{\prime}}(y)}{x-y}dy=p_{\theta^{\prime}}^{\prime}(x).
\end{align*}
Therefore, $h$ and $p_{\theta^{\prime}}$ differ only by an additive
constant, which may be absorbed into the normalization $\phi(n)$.
Thus $h(A)=p(A)+\sum_{i=1}^m\theta^{\prime}_i\frac{\partial}{\partial\theta_i}p_{\theta,n}(A)$,
and we conclude that all curves in $S$ throught $\theta$ are given by a linear combination of $\frac{\partial}{\partial\theta_i}p_{\theta,n}(A)$, $i=1,\ldots,m$.
\end{proof}

~

\subsection{The Fisher Information Metric\label{sec:The-Metric}}

~

A key feature of the metric defined by \cite{AmNa} is that it
satisfies the equation
\begin{align}\label{eq:crucial-equation}
g_{ij}(\theta)=\frac{\partial^{2}}{\partial\theta_{i}\partial\theta_{j}}\psi(\theta).
\end{align}

We take this equation as the starting point for our definitions, and
we calculate
\begin{align}\label{eq:derivative-of-pressure}
\frac{\partial}{\partial\theta_{j}}\psi(\theta,n)=\left(\frac{1}{n^{2}}\frac{\int-n\textrm{Tr}\left(F_{j}\right)\exp\left(-n\Tr\left(p(A)+\sum_{k=1}^m\theta_{k}F_{k}(A)\right)\right)dA}{\int\exp\left(-n\Tr\left(p(A)+\sum_{k=1}^m\theta_{k}F_{k}(A)\right)\right)dA}\right)=-\frac{1}{n}\int\textrm{Tr}\left(F_{j}\right)d\mu_{\theta,n}(A),
\end{align}
and
\begin{align}\label{eq:metric-equation}
\frac{\partial^2}{\partial\theta_i\partial\theta_j}\psi(\theta,n)=\int\Tr\Big(F_j\Big)\Tr\Big(F_i+\frac{\partial}{\partial\theta_i}\psi(\theta,n)\Big)d\mu_{\theta,n}(A).
\end{align}

~

We define an inner-product and a corresponding metric and check
that it satisfies \eqref{eq:metric-equation}.
\begin{defn}
\begin{align*}
\left\langle{f(A),g(A)}\right\rangle_{\theta,n}=\int\textrm{Tr}\left(f(A)\right)\textrm{Tr}\left(g(A)\right)d\mu_{\theta,n}(A)
\end{align*}
\begin{align}\label{eq:metric}
g_{ij}(\theta,n)=\left\langle
\frac{\partial}{\partial\theta_{i}}p_{\theta,n},\frac{\partial}{\partial\theta_{j}}p_{\theta,n}\right\rangle_{\theta,n}=\int\textrm{Tr}\left(\frac{\partial}{\partial\theta_{i}}p_{\theta,n}(A)\right)\textrm{Tr}\left(\frac{\partial}{\partial\theta_{j}}p_{\theta,n}(A)\right)d\mu_{\theta,n}(A).
\end{align}
\end{defn}

~

\begin{prop}
$g_{ij}(\theta,n)=\frac{\partial^{2}}{\partial\theta_{i}\partial\theta_{j}}\psi(\theta,n)$.
\end{prop}

\begin{proof}
We calculate:
\begin{align*}
g_{ij}\left(\theta,n\right)&=\int\textrm{Tr}\left(F_{i}+\frac{\partial}{\partial\theta_{i}}\psi(\theta,n)\right)\textrm{Tr}\left(F_{j}+\frac{\partial}{\partial\theta_{j}}\psi(\theta,n)\right)d\mu_{\theta,n}(A)\\
&=\int\textrm{Tr}\left(F_{i}-\frac{1}{n}\int\textrm{Tr}(F_{i})d\mu_{\theta,n}(A)\right)\textrm{Tr}\left(F_{j}-\frac{1}{n}\int\textrm{Tr}(F_{j})d\mu_{\theta,n}(A)\right)d\mu_{\theta,n}(A),
\end{align*}
by equation \eqref{eq:derivative-of-pressure}.

~

Now $\Tr\left(\frac{1}{n}\int\textrm{Tr}(F_{i})d\mu_{\theta,n}(A)\right)=\int\textrm{Tr}(F_{i})d\mu_{\theta,n}(A)$, so we have
\begin{align*}
g_{ij}\left(\theta,n\right)&=\int\textrm{Tr}\left(F_{i}\right)\textrm{Tr}\left(F_{j}\right)d\mu_{\theta,n}(A)-\int\textrm{Tr}\left(F_{i}\right)d\mu_{\theta,n}(A)\int\textrm{Tr}\left(F_{j}\right)d\mu_{\theta,n}(A)\\
&=\int\textrm{Tr}\left(F_{i}\right)\textrm{Tr}\left(F_{j}+\frac{\partial}{\partial\theta_{j}}\psi(\theta,n)\right)d\mu_{\theta,n}(A)\\
&=\frac{\partial^{2}}{\partial\theta_{i}\partial\theta_{j}}\psi(\theta,n),
\end{align*}
using equations \eqref{eq:derivative-of-pressure} and \eqref{eq:metric-equation}.
\end{proof}

~

\begin{rem}
In connection to Voiculescu's free probability theory, it seems
natural to define the inner-product on the tangent space
as
\begin{align}\label{eq:badmetric}
\left\langle{f,g}\right\rangle_{\theta}=\frac{1}{n}\int\Tr\left(f(A)g(A)\right)d\mu_{\theta,n}(A).
\end{align}
However, the metric \eqref{eq:badmetric} does not satisfy equation
\eqref{eq:crucial-equation}, which is crucial in order to calculate
the Legendre transform in Section \ref{sec:Legendre-Transform}.

~

Also, if we restrict our attention to Gaussian Unitary Ensemble
perturbations (which give rise to a semicircular perturbation as
$n\rightarrow\infty$), the metric \eqref{eq:metric}
coincides with metric \eqref{eq:badmetric} (see Section
\ref{sec:Conjugate-Variable}).
\end{rem}

~

\subsection{The $(\alpha)$-Connections\label{sec:The-alpha-Connection}}

~

To calculate the Legendre transform we need a pair of dual affine
connections on the manifold (see \cite{BiGo}). In fact we define a
family of pairs of dual connections with curvature parameter
$\alpha\in\left[-1,1\right]$, and we denote the connection
coefficients by $\Gamma_{ijk}^{(\alpha),n}\left(\theta\right)$. To
compute the Legendre transform, the $\left(\alpha\right)$-connection
must be dual to the $\left(-\alpha\right)$-connection, i.e.\begin{align*}
\frac{\partial}{\partial\theta_{k}}g_{ij}\left(\theta,n\right)=\Gamma_{kij}^{\left(\alpha\right),n}\left(\theta\right)+\Gamma_{kji}^{\left(-\alpha\right),n}\left(\theta\right)\textrm{.}\end{align*}

We take this as our starting point for the definition of $\Gamma$,
and we calculate:
\begin{align*}
&\frac{\partial}{\partial\theta_{k}}g_{ij}\left(\theta,n\right)=\frac{\partial^{3}}{\partial\theta_{i}\partial\theta_{j}\partial\theta_{k}}\psi\left(\theta,n\right)=\\
&\int\Tr\left(\frac{\partial^{2}}{\partial\theta_{k}\partial\theta_{i}}p_{\theta,n}(A)\right)\textrm{Tr}\left(\frac{\partial}{\partial\theta_{j}}p_{\theta,n}(A)\right)d\mu_{\theta,n}(A)+\int\Tr\left(\frac{\partial}{\partial\theta_{i}}p_{\theta,n}(A)\right)\textrm{Tr}\left(\frac{\partial^{2}}{\partial\theta_{k}\partial\theta_{j}}p_{\theta,n}(A)\right)d\mu_{\theta,n}(A)\\
&-n\int\Tr\left(\frac{\partial}{\partial\theta_{i}}p_{\theta,n}(A)\right)\textrm{Tr}\left(\frac{\partial}{\partial\theta_{j}}p_{\theta,n}(A)\right)\textrm{Tr}\left(\frac{\partial}{\partial\theta_{k}}p_{\theta,n}(A)\right)d\mu_{\theta,n}(A).
\end{align*}

This leads us to the following definition

~

\begin{defn}
\begin{align*}
\Gamma_{ijk}^{\left(\alpha\right),n}\left(\theta\right)=&\int\textrm{Tr}\left(\frac{\partial^{2}}{\partial\theta_{k}\partial\theta_{i}}p_{\theta,n}(A)\right)\textrm{Tr}\left(\frac{\partial}{\partial\theta_{j}}p_{\theta,n}(A)\right)d\mu_{\theta,n}(A)\\
&-\frac{1-\alpha}{2}\cdot{n}\int\textrm{Tr}\left(\frac{\partial}{\partial\theta_{i}}p_{\theta,n}(A)\right)\textrm{Tr}\left(\frac{\partial}{\partial\theta_{j}}p_{\theta,n}(A)\right)\textrm{Tr}\left(\frac{\partial}{\partial\theta_{k}}p_{\theta,n}(A)\right)d\mu_{\theta,n}(A).
\end{align*}
\end{defn}

~

Notice that the connection coefficients depend on the choice of
coordinate system, and this gives the notion of flatness:
\begin{defn}
A coordinate system $\left\{\zeta_{i}\right\}$ is
$\left(\alpha\right)$-flat if
$\Gamma_{ijk}^{\left(\alpha\right),n}\left(\zeta\right)=0$.
That is, the vector fields $X_{i}=\frac{\partial}{\partial\zeta_{i}}$ are parallel with respect to the $\left(\alpha\right)$-connection.
\end{defn}

\begin{prop}
An exponential family $S$ is $(1)$-flat.
\end{prop}
\begin{proof}
Calculate:
\begin{align*}
\Gamma_{ijk}^{\left(1\right),n}\left(\theta\right)&=\int\textrm{Tr}\left(\frac{\partial^{2}}{\partial\theta_{k}\partial\theta_{i}}p_{\theta,n}(A)\right)\textrm{Tr}\left(\frac{\partial}{\partial\theta_{j}}p_{\theta,n}(A)\right)d\mu_{\theta,n}(A)\\
&=\int\textrm{Tr}\left(\frac{\partial}{\partial\theta_{k}}\left(F_{i}(A)-\frac{\partial}{\partial\theta_{i}}\psi\left(\theta,n\right)\right)\right)\textrm{Tr}\left(F_{j}(A)-\frac{\partial}{\partial\theta_{j}}\psi\left(\theta,n\right)\right)d\mu_{\theta,n}(A)\\
&=\int\textrm{Tr}\left(-\frac{\partial^{2}}{\partial\theta_{k}\partial\theta_{i}}\psi\left(\theta,n\right)\right)\textrm{Tr}\left(F_{j}(A)-\frac{\partial}{\partial\theta_{j}}\psi\left(\theta,n\right)\right)d\mu_{\theta,n}(A)\\
&=-n\frac{\partial^{2}}{\partial\theta_{k}\partial\theta_{i}}\psi\left(\theta,n\right)\cdot\int\textrm{Tr}\left(F_{j}(A)-\frac{1}{n}\int\textrm{Tr}\left(F_{j}(A)\right)d\mu_{\theta,n}(A)\right)d\mu_{\theta,n}(A)\\
&=-n\frac{\partial^{2}}{\partial\theta_{k}\partial\theta_{i}}\psi\left(\theta,n\right)\cdot\left(\int\textrm{Tr}\left(F_{j}(A)\right)d\mu_{\theta,n}(A)-\int\textrm{Tr}\left(F_{j}(A)\right)d\mu_{\theta,n}(A)\right)=0.
\end{align*}
\end{proof}
\begin{prop}\label{pro: dual-connections}
The $\left(\alpha\right)$- and $\left(-\alpha\right)$-connections
are mutually dual with respect to $g\left(\theta,n\right)$.
\end{prop}
\begin{proof}
This follows from the calculation:
\begin{align*}
\frac{\partial}{\partial\theta_{k}}g_{ij}\left(\theta,n\right)=&\int\textrm{Tr}\left(\frac{\partial^{2}}{\partial\theta_{k}\partial\theta_{i}}p_{\theta,n}(A)\right)\textrm{Tr}\left(\frac{\partial}{\partial\theta_{j}}p_{\theta,n}(A)\right)d\mu_{\theta,n}(A)\\
&-\int\textrm{Tr}\left(\frac{\partial}{\partial\theta_{i}}p_{\theta,n}(A)\right)\textrm{Tr}\left(\frac{\partial^{2}}{\partial\theta_{k}\partial\theta_{j}}p_{\theta,n}(A)\right)d\mu_{\theta,n}(A)\\
&-\frac{1-\alpha}{2}\cdot{n}\cdot\int\textrm{Tr}\left(\frac{\partial}{\partial\theta_{i}}p_{\theta,n}(A)\right)\textrm{Tr}\left(\frac{\partial}{\partial\theta_{j}}p_{\theta,n}(A)\right)\textrm{Tr}\left(\frac{\partial}{\partial\theta_{k}}p_{\theta,n}(A)\right)d\mu_{\theta,n}(A)\\
&+\frac{1+\alpha}{2}\cdot{n}\cdot\int\textrm{Tr}\left(\frac{\partial}{\partial\theta_{i}}p_{\theta,n}(A)\right)\textrm{Tr}\left(\frac{\partial}{\partial\theta_{j}}p_{\theta,n}(A)\right)\textrm{Tr}\left(\frac{\partial}{\partial\theta_{k}}p_{\theta,n}(A)\right)d\mu_{\theta,n}(A)\\
=&\Gamma_{ijk}^{\left(\alpha\right),n}\left(\theta\right)+\Gamma_{ikj}^{\left(-\alpha\right),n}\left(\theta\right).
\end{align*}
\end{proof}

~

\begin{cor}
$\frac{\partial}{\partial\theta_{k}}g_{ij}\left(\theta,n\right) =
\Gamma_{ijk}^{\left(0\right),n}\left(\theta\right)+\Gamma_{ikj}^{\left(0\right),n}\left(\theta\right)$,
i.e. the $\left(0\right)$-connection is the Levi-Cevita (metric) connection.
\end{cor}
~

\subsection{Several Independent Matrices\label{sec:Several-Independent-Matrices}}

~

 Our discussion started with a single random matrix model to
familiarize the reader with the geometric notions and calculations,
but in fact it extends to several independent matrices as follows. Start with independent random matrices $A_{1},\ldots,A_{k}$ with $A_{r}$ distributed according to
\begin{align*}
\exp\left(-n\textrm{Tr}\left(p\left(r\right)\left(A\right)+\sum_{i=1}^{m}\theta_{i}\left(r\right)F_{i}\left(r\right)\left(A\right)+\psi\left(r\right)\left(\theta\left(r\right),n\right)\right)\right)dA.
\end{align*}
Then $\left(A_{1},\ldots,A_{k}\right)$ is distributed according to
\begin{align*}
\exp\left(-n\textrm{Tr}\left(\left(\sum_{r=1}^{k}p\left(r\right)\left(A_{r}\right)\right)+\left(\sum_{r=1}^{k}\sum_{i=1}^{m}\theta_{i}\left(r\right)F_{i}\left(r\right)\left(A_{r}\right)\right)+\left(\sum_{r=1}^{k}\psi\left(r\right)\left(\theta\left(r\right),n\right)\right)\right)\right)dA_{1}\ldots
dA_{k}.
\end{align*}

~

\begin{notation}
Write $P\left(A_{1},\ldots,A_{k}\right)$ for $\sum_{r=1}^{k}p\left(r\right)\left(A_{r}\right)$, $\theta$ for\\ $\left(\theta_{1}\left(1\right),\ldots,\theta_{m}\left(1\right),\theta_{1}\left(2\right),\ldots,\theta_{m}\left(2\right),\ldots\ldots,\theta_{m}\left(k\right)\right)$, and $P_{\theta,n}\left(A_{1},\ldots,A_{k}\right)$ for\\ $P\left(A_{1},\ldots,A_{k}\right)+\left(\sum_{r=1}^{k}\sum_{i=1}^{m}\theta_{i}\left(r\right)F_{i}\left(r\right)\left(A_{r}\right)\right)+\Psi\left(\theta,n\right)$,\\ with $\Psi\left(\theta,n\right)=\frac{1}{n^{2}}\log\int\exp\left(-n\textrm{Tr}\left(P\left(A_{1},\ldots,A_{k}\right)+\left(\sum_{r=1}^{k}\sum_{i=1}^{m}\theta_{i}\left(r\right)F_{i}\left(r\right)\left(A_{r}\right)\right)\right)\right)dA_{1}\ldots{dA_k}$. Also write $d\mu_{\theta(r),n}(A_{r})$ for $\exp\left(-n\textrm{Tr}\left(p\left(r\right)\left(A_{r}\right)+\sum_{i=1}^{m}\theta_{i}\left(r\right)F_{i}\left(r\right)\left(A_{r}\right)\right)\right)dA_{r}$, and $d\tilde{\mu}_{\theta,n}(A_{1},\ldots,A_{k})$ for $d\mu_{\theta(1),n}(A_{1})\cdots{d}\mu_{\theta(r),n}(A_{k})$.
\end{notation}

~

Notice, by recentering and rescaling we may assume that $\frac{1}{n}\int\textrm{Tr}\left(A_{r}\right)d\mu_{\theta(r),n}(A_{r})=0$ and\\ $\frac{1}{n}\int\textrm{Tr}\left(A_{r}^{2}\right)d\mu_{\theta(r),n}(A_{r})=1$ for $r=1,\ldots,k$.  Since $A_1,\ldots,A_k$ are independent,
\cite[Theorem 4.4.1]{VDN} shows that
$\left(A_{1},\ldots,A_{k}\right)$ converges as $n\rightarrow\infty$
in the sense that $\exists\tau$ a tracial state on
$\mathbb{C}\left\langle x_1,\ldots,x_k\right\rangle$ such that
$\forall w\in\mathbb{C}\left\langle x_1,\ldots,x_k\right\rangle$
\begin{align*}
\frac{1}{n}\int\Tr\left(w\left(A_{1},\ldots,A_{k}\right)\right)d\tilde{\mu}_{\theta,n}\left(A_{1},\ldots,A_{k}\right)
\rightarrow\tau\left(w\right)
\end{align*}

Now we must identify the potential. In the classical case, when one
considers several independent random variables $X_{1},\ldots,X_{k}$
instead of a single variable, with $X_i$ distributed according
to $\exp\left(p_{i}(x)\right)dx$, the potential is determined by the
equation
\begin{align*}
\frac{d}{dx_{r}}p\left(x_{1},\ldots,x_{k}\right)=d_{x_{r}}^{*}(1),
\end{align*}
for $r=1,\ldots,k$, with
\begin{align*}
d_{x_{r}}:L^{2}\left(\exp\left(\sum_{r=1}^{k}p_{r}(x)\right)dx_1{\ldots}dx_k\right)&\rightarrow{L}^2\left(\exp\left(\sum_{r=1}^{k}p_{r}(x)\right)dx_1{\ldots}dx_k\right),\\
g(x_{1},\ldots,x_{k})&\mapsto\frac{d}{dx_{r}}g(x_1,\ldots,x_k)
\end{align*}
the densely defined partial differentiation operator with
$\textrm{domain}(d_{x_{r}})=\textrm{polynomials}$.

In the multi-matrix case, the equation for $h_{\theta}$ to be a
potential in the limit $n\rightarrow\infty$ becomes
\begin{align*}
\mathcal{D}_{x_{r}}h_{\theta}\left(x_{1},\ldots,x_{k}\right)=\partial_{x_{r}}^{*}\left(1\otimes1\right),
\end{align*}
with $\partial_{x_{r}}$ and $\mathcal{D}_{x_{r}}$ defined in
\cite{Voi2} as follows. $\partial_{x_{r}}:L^{2}\left(\mathbb{C}\left\langle{x_1,\ldots,x_k}\right\rangle,\tau\right)\rightarrow{L}^2\left(\mathbb{C}\left\langle{x_1,\ldots,x_k}\right\rangle,\tau\right)\otimes{L}^2\left(\mathbb{C}\left\langle{x_1,\ldots,x_k}\right\rangle,\tau\right)$ is defined by $\partial_{x_{r}}(x_{s})=\delta_{rs}\cdot1\otimes1$ and $\partial_{x_{r}}(1)=0$. $\mathcal{D}_{x_r}:L^2(\mathbb{C}\langle{x_1,\ldots,x_k}\rangle,\tau)\rightarrow{L}^2(\mathbb{C}\langle{x_1,\ldots,x_k}\rangle,\tau)$ is defined by $\mathcal{D}_{x_r}=\sigma\circ\partial_{x_{r}}$ where $\sigma(x\otimes{y})=yx$.

~

\cite[Prop 3.6]{Voi2} shows that
$\partial_{x_r}^*\left(1\otimes1\right)\in{L}^2\left(\mathbb{C}\langle{x_r}\rangle,\tau\right)$, so we may apply the discussion from Section \ref{sec:The-Tangent-Space} to get
\begin{align*}
\partial_{x_{r}}^{*}\left(1\otimes1\right)=\frac{d}{dx_r}\left(p\left(r\right)\left(x_{r}\right)+\sum_{i=1}^{m}\theta_{i}\left(r\right)F_{i}\left(r\right)\left(x_{r}\right)\right).
\end{align*}

~

Thus, the condition for a potential of several independent matrices
becomes
\begin{align*}
\mathcal{D}_{x_{r}}h\left(x_{1},\ldots,x_{k}\right)=\frac{d}{dx_r}\left(p\left(r\right)\left(x_{r}\right)+\sum_{i=1}^{m}\theta_{i}\left(r\right)F_{i}\left(r\right)\left(x_{r}\right)\right).
\end{align*}

~

This has a solution
\begin{align*} h\left(x_{1},\ldots,x_{k}\right) =
P_{\theta,n}\left(x_{1},\ldots,x_{k}\right)\textrm{,}\end{align*} which leads
to the definition

~

\begin{defn}
The potential of $d\tilde{\mu}_{\theta,n}\left(A_{1},\ldots,A_{k}\right)$
is $P_{\theta,n}\left(A_{1},\ldots,A_{k}\right)$.
\end{defn}

~

\begin{defn}
The tangent space to this model is $T_{\theta}S=\mathrm{span}\left\{ \frac{\partial}{\partial\theta_{i}}P_{\theta,n}\left(A_{1},\ldots,A_{k}\right)\right\} _{i=1}^{m}\textrm{.}$
\end{defn}

~

Next we define an inner-product and a corresponding metric that
satisfies equation (\ref{eq:crucial-equation}).

~

\begin{defn}
\begin{align*}
\left\langle f,h \right\rangle _{\theta}=
-\int\textrm{Tr}\left(f\left(A_{1},\ldots,A_{k}\right)\right)\textrm{Tr}\left(h\left(A_{1},\ldots,A_{k}\right)\right)d\tilde{\mu}_{\theta,n}\left(A_{1},\ldots,A_{k}\right)\textrm{.}\end{align*}
\end{defn}

~

\begin{defn}
\begin{align*}
G_{i(r)j(s)}\left(\theta,n\right)=\int\textrm{Tr}\left(\frac{\partial}{\partial\theta(r)_{i}}P_{\theta,n}\right)\textrm{Tr}\left(\frac{\partial}{\partial\theta(s)}P_{\theta,n}\right)d\tilde{\mu}_{\theta,n}\left(A_{1},\ldots,A_{k}\right)\textrm{.}\end{align*}
\end{defn}

It is a straightforward calculation that

\begin{prop}
$G_{i(r)j(s)}\left(\theta,n\right)=\frac{\partial^{2}}{\partial\theta(r)_{i}\partial\theta(s)_{j}}\Psi\left(\theta,n\right)$.
\end{prop}

~

\begin{rem}\label{rem:independent-matrices}
Notice that for $r\neq s$ $A_{r}$ is independent of
$A_{s}$, so $G_{i(r)j(s)}\left(\theta,n\right)=0$. Therefore, if
$A_{1},\ldots,A_{k}$ are identically distributed we have
$p_{1}=p_{2}=\ldots=p_{k}$ and
$\theta\left(1\right)=\theta\left(2\right)=\ldots=\theta\left(k\right)$,
so
\begin{align*}
\Psi\left(\theta,n\right)&=\frac{1}{n^{2}}\log\int\exp\left(-n\textrm{Tr}\left(\sum_{r=1}^{k}\left(p_{r}\left(A_{r}\right)+\sum_{i=1}^{m}\theta_{i}\left(r\right)F_{i}\left(A_{r}\right)\right)\right)\right)dA_{1}\ldots{d}A_k\\
&=\frac{1}{n^{2}}\sum_{r=1}^{k}\log\int\exp\left(-n\textrm{Tr}\left(p_{r}\left(A_{r}\right)+\sum_{i=1}^{m}\theta_{i}\left(r\right)F_{i}\left(A_{r}\right)\right)\right)dA_{r}.
\end{align*}
Thus, as in the classical case

\begin{align*}
G\left(\theta,n\right)&=
\underbrace{\left(
    \begin{array}{ccc}
    g(\theta(1),n) &  & 0\\
    & \ddots & \\
    0 &  & g(\theta(1),n)
    \end{array}
\right)}_k.
\end{align*}
\end{rem}

~

We also define the $(\alpha)$-connections
\begin{defn}
\begin{align*}
\Gamma_{i(r)j(s)k(t)}^{(\alpha),n}\left(\theta\right)=&\int\textrm{Tr}\left(\frac{\partial^{2}}{\partial\theta(t)_{k}\partial\theta(r)_{i}}P_{\theta,n}\right)\textrm{Tr}\left(\frac{\partial}{\partial\theta(s)_{j}}P_{\theta,n}\right)d\tilde{\mu}_{\theta,n}(A_{1},\ldots,A_{k})\\
&-\frac{1-\alpha}{2}\cdot{n}\int\textrm{Tr}\left(\frac{\partial}{\partial\theta(r)_{i}}P_{\theta,n}\right)\textrm{Tr}\left(\frac{\partial}{\partial\theta(s)_{j}}P_{\theta,n}\right)\textrm{Tr}\left(\frac{\partial}{\partial\theta(t)_{k}}P_{\theta,n}\right)d\tilde{\mu}_{\theta,n}(A_{1},\ldots,A_{k}).
\end{align*}
\end{defn}
Using the argument in Remark \ref{rem:independent-matrices}, for identically distributed $A_1,\ldots,A_k$ we have $\Gamma_{i(r)j(s)k(t)}^{(\alpha),n}\left(\theta\right)=0$
if $r\neq s$, $s\neq t$, or $r\neq t$, and $\Gamma_{i(r)j(s)k(t)}^{(\alpha),n}\left(\theta\right)=\Gamma_{ijk}^{(\alpha),n}(\theta(1))$ otherwise.

~

\begin{defn}
Given independent random matrices $A_1,\ldots,A_k$ with
distributions of the form \eqref{eq:exponential}, the
\emph{Information Manifold} associated to $\left(A_1,\ldots,A_k\right)$ is the geometric structure $S=\left(M,g,\Gamma,\Gamma^*\right)$ described in this section.
\end{defn}
Combining our observations, we obtain the following theorem,

\begin{thm}\label{thm:manifold-of-independent-matrices} Let $A_1,\ldots,A_k$ be independent random matrices
with distribution functions of the form (\ref{eq:exponential}), let
$S_i$ be the information manifold associated to $A_i$, and let
$\tilde{S}$ be the information manifold associated to
$\left(A_1,\ldots,A_k\right)$, then $\tilde{S}=S_1\oplus\ldots\oplus{S}_k$.
\end{thm}

~

\subsection{Legendre Transform of Pressure\label{sec:Legendre-Transform}}

~

 In this section we calculate the Legendre transform of the
pressure $\psi\left(\theta,n\right)$. Since the notation for several
independent matrices is cumbersome, we calculate with a single
matrix and by theorem (\ref{thm:manifold-of-independent-matrices})
our calculations extend to several independent matrices.

\cite[Section 3.3-3.5]{AmNa} define the Legendre transform for a
smooth real-valued function on a Reimannian manifold with a pair of
dual connections; this includes our construction, so their
discussion applies in our case.

Following \cite[Section 3.3-3.5]{AmNa}, define a new coordinate
system
\begin{align*}
\eta_{i}=\frac{1}{n}\int\textrm{Tr}\left(F_{i}(A)\right)d\mu_{\theta,n}(A)\textrm{.}\end{align*}
According to equation (\ref{eq:derivative-of-pressure}),
\begin{align*}\frac{\partial}{\partial\theta_{i}}\psi(\theta,n)=\frac{1}{n}\int\textrm{Tr}\left(F_{i}(A)\right)d\mu_{\theta,n}(A)=\eta_{i}\textrm{,}\end{align*}
so we have
\begin{align*}\frac{\partial\eta_{j}}{\partial\theta_{i}}=\frac{\partial^{2}}{\partial\theta_{i}\partial\theta_{j}}\psi=g_{ij}\left(\theta,n\right)\textrm{,}\end{align*}
and similarly\begin{align*}
\frac{\partial\theta_{j}}{\partial\eta_{i}}=\left(g^{-1}\left(\theta,n\right)\right)_{ij}\textrm{.}\end{align*}
Therefore, $\left\{ \theta_{i}\right\} $ and $\left\{
\eta_{i}\right\} $ are coordinate systems which are mutually dual
with respect to $g_{ij}\left(\theta,n\right)$ \cite[Section
3.5]{AmNa}.

~

\begin{prop}
$\left\{ \eta_{i}\right\} $ is $(-1)$-flat.
\end{prop}
\begin{proof}
The $\left\{ \theta_{i}\right\} $ coordinate system is
$\left(1\right)$-flat. By proposition (\ref{pro: dual-connections}),
the dual coordinate system is $\left(-1\right)$-flat, so $\left\{
\eta_{i}\right\} $ is $(-1)$-flat
\end{proof}

~

\begin{defn}
The \emph{Legendre Transform} of $\psi\left(\theta,n\right)$ is
\begin{align*}
\varphi\left(\theta,n\right)=\sup_{\theta^{\prime}\in\Theta}\left\{
\frac{1}{n}\int\textrm{Tr}\left(p(A)+\sum_{i=1}^{m}\theta_{i}^{\prime}F_{i}\left(A\right)+\psi\left(\theta^{\prime},n\right)\right)d\mu_{\theta,n}(A)-
\frac{1}{n}\int\textrm{Tr}\left(p(A)\right)d\mu_{\theta,n}(A)\right\}
\textrm{.}\end{align*}
\end{defn}

~

\begin{prop}\label{pro:legendre-transform}
\begin{align*}
\varphi\left(\theta,n\right)=\frac{1}{n}\int\textrm{Tr}\left(p_{\theta,n}(A)\right)d\mu_{\theta,n}(A)-\frac{1}{n}\int\textrm{Tr}\left(p(A)\right)d\mu_{\theta,n}(A)\textrm{.}\end{align*}

\end{prop}
\begin{proof}
\cite[Section 3.5]{AmNa} show that
\begin{align*}
\varphi(\theta,n) =
\sum_{i=1}^{m}\theta_{i}\eta_{i}(\theta)+\psi(\theta,n)=\frac{1}{n}\int\textrm{Tr}\left(p_{\theta,n}(A)\right)d\mu_{\theta,n}(A)-\frac{1}{n}\int\textrm{Tr}\left(p(A)\right)d\mu_{\theta,n}(A).\end{align*}
\end{proof}

~

Following the discussion in \cite[Section 3.5]{AmNa}, the first
summand in proposition \eqref{pro:legendre-transform} is the
analogue of entropy,

\begin{defn}
\begin{align*}
H(d\mu_{\theta,n})=-\frac{1}{n}\int\Tr\left(p_{\theta,n}(A)\right)d\mu_{\theta,n}(A).
\end{align*}
\end{defn}

Also according to the discussion in \cite[Section 3.5]{AmNa},

\begin{cor}\label{cor:convex}
Both $H(d\mu_{\theta,n})$ and $\psi(\theta,n)$ are
convex in $\theta$.
\end{cor}

~

\subsection{Calculations on Well-Known Models.\label{sec:Calculations}}

~

 In this section we present some calculations on random matrix
models that appear in applications.

\subsubsection{Gaussian Unitary Ensemble (GUE)}

~

~

Classically the most ubiquitous statistical model is the Gaussian
family,
\begin{align*}\exp\left(-\frac{(x-\mu)^{2}}{2\sigma^{2}}-\frac{1}{2}\log(2\pi\sigma^{2})\right),\end{align*}
which is an exponential family by setting $\theta_{1}=-\frac{\mu}{\sigma^{2}}$, $\theta_{2}=\frac{1}{2\sigma^2}$, $\psi(\theta_{1},\theta_{2})=\frac{\left(\theta_{1}\right)^{2}}{2\theta_{2}}+\frac{1}{2}\log\frac{\pi}{\theta_{2}}$,
and writing
\begin{align*}
\exp\left(-\frac{(x-\mu)^{2}}{2\sigma^{2}}-\frac{1}{2}\log(2\pi\sigma^{2})\right)=\exp\left(-\left(\theta_{2}x^{2}+\theta_{1}x+\psi(\theta_{1},\theta_{2})\right)\right).
\end{align*}

The corresponding random matrix model is the Gaussian Unitary
Ensemble (GUE), with distribution function
\begin{align*}
\exp\left(-n\textrm{Tr}\left(\frac{\left(A-\mu\right)^{2}}{2\sigma^{2}}+\frac{1}{2}\log\frac{2\pi\sigma^{2}}{n}\right)\right),
\end{align*}
which is an exponential family under the same coordinates
$\tilde{\theta}_{1}=-\frac{\mu}{\sigma^{2}}$,
$\tilde{\theta}_{2}=\frac{1}{2\sigma^{2}}$, and
$\tilde{\psi}(\tilde{\theta}_{1},\tilde{\theta}_{2},n)=\left.\left(\tilde{\theta}_{1}\right)^2\right\slash
2\tilde{\theta}_2+\frac{1}{2}\log\left(\pi\slash
n\tilde{\theta}_2\right)$.

We see that
\begin{align*}
\frac{\partial^{2}}{\partial\theta_{i}\partial\theta_{j}}\psi(\theta_{1},\theta_{2})=\frac{\partial^{2}}{\partial\tilde{\theta}_{i}\partial\tilde{\theta}_{j}}\tilde{\psi}(\tilde{\theta}_{1},\tilde{\theta}_{2},n),
\end{align*}
so the Fisher information metric of the GUE model in $\left(\tilde{\theta}_{1},\tilde{\theta}_{2}\right)$-coordinates
is the same as the classical Fisher information metric of the Gaussian
model in $\left(\theta_{1},\theta_{2}\right)$-coordinates. Using the change of coordinate rule $g_{kl}(\xi)=\sum_{i,j=1}^{m}g_{ij}(\rho)\left(\frac{\partial\rho^{i}}{\partial\xi_{k}}\right)\left(\frac{\partial\rho^{j}}{\partial\xi_{l}}\right)$, the equality $\frac{\partial\theta_{i}}{\partial\mu}=\frac{\partial\tilde{\theta}_{i}}{\partial\mu}$, and the equality $\frac{\partial\theta_{i}}{\partial\sigma}=\frac{\partial\tilde{\theta}_{i}}{\partial\sigma}$ for $i=1,2$, we see that the Fisher information metric for the GUE model in the
$(\mu,\sigma)$-coordinates is the same as the Fisher information
metric for the Gaussian model in the $(\mu,\sigma)$-coordinates:
\begin{align*}
g_{ij}(\mu,\sigma,n)=\left(
    \begin{array}{cc}
    \frac{1}{\sigma^2} & 0\\
    0 & \frac{2}{\sigma^{2}}
    \end{array}
\right).
\end{align*}

Taking the limit $n\rightarrow\infty$, the semicircle
and Gaussian distributions have the same Fisher information metric.

\subsubsection{Laguerre Unitary Ensemble (LUE)}

~

~

Another well-known model in random matrix theory is the Laguerre
Unitary Ensemble (LUE) also known as a Wishart matrix. This is
defined in \cite{FoWi} as the random matrix $A=X^*X$ with $X\sim\exp\left(-n\textrm{Tr}\left(X^{2}\right)\right)dX$ (we
take $n=N$ in the definition of \cite{FoWi} so their notation matches ours). \cite{FoWi} show that the eigenvalues of $A$ are distributed according to
\begin{align*}
Z_n^{-1}\cdot\prod_{1\leq{i}<j\leq{n}}\left(\lambda_{i}-\lambda_{j}\right)\prod_{k=1}^{n}\exp\left(\lambda_{k}\right)\cdot\chi_{\lambda_{k}>0}d\lambda_{k},
\end{align*}
with $Z_{n}$ the normalization constant. Therefore, $A$ is
distributed according to $Z_{n}^{-1}\cdot\exp\left(-n\textrm{Tr}\left(A\right)\right)dA$ on $\left\{ \left.A\in M_{n}^{SA}(\mathbb{C})\right|A>0\right\}$ as an
orthogonally invariant model. We will not discuss orthogonally-invariant models in general, but our definitions make sense verbatim in this case.

Instead of starting with a standard LUE, we may parameterize its
variance, and rescale it for convergence as $n\rightarrow\infty$ to
obtain
\begin{align*}
Z_n^{-1}(t)\cdot\exp\left(-n\textrm{Tr}\left(\frac{A}{t}\right)\right)\textrm{ on }\left\{\left.A\in M_{n}^{SA}(\mathbb{C})\right|A>0\right\}
\end{align*}
with
\begin{align*}
Z_n(t)=\mathop{\int}_{\left\{\left.A\in{M}_n^{SA}(\mathbb{C})\right|A>0\right\}}\exp\left(-n\textrm{Tr}\left(\frac{A}{t}\right)\right)dA=\int_{M_{n}^{SA}(\mathbb{C})}\exp\left(-n\textrm{Tr}\left(\frac{A^{2}}{t}\right)\right)dA=\frac{1}{2}\log\frac{\pi{t}}{n}.
\end{align*}
This is an exponential family by setting
$\theta_{1}=1\slash{t}$, $\psi(\theta_{1},n)=\frac{1}{2}\log\left(\pi\slash{n}\theta_{1}\right)$, and writing
\begin{align*}
Z_{n}^{-1}(t)\cdot\exp\left(-n\textrm{Tr}\left(\frac{A}{t}\right)\right)=\exp\left(-n\textrm{Tr}\left(\theta_{1}A+\psi(\theta_{1},n)\right)\right).
\end{align*}
Now $\frac{\partial^2}{\partial\theta_1\partial\theta_1}\psi(\theta_1,n)=\frac{1}{2\left(\theta_1\right)^2}$, so the Fisher information metric of the LUE model is
\begin{align*}
g_{11}(t)=\frac{1}{2t^{2}}.
\end{align*}

~

~

\section{The $n\rightarrow\infty$ case}\label{sec:The-n-goes-to-infinity-case}

\subsection{Convergence}

~

In this section we verify that the tangent space, Fisher information
metric, $\left(\alpha\right)$-connections, pressure, and entropy
converge as $n\rightarrow\infty$. In fact we find that the entropy
and pressure converge to the free entropy and free pressure, and the
Fisher information metric of the semicircular perturbation model
coincides with Voiculescu's Fisher information measure.

~

First we note that the tangent space, regarded as a vector space of
random variables, converges in moments:

\begin{prop}
\begin{align*}
\frac{1}{n}\int\textrm{Tr}\left(\frac{\partial}{\partial\theta_{i_{1}}}p_{\theta,n}(A)\cdots\frac{\partial}{\partial\theta_{i_{k}}}p_{\theta,n}(A)\right)d\mu_{\theta,n}(A)\rightarrow
\int\frac{\partial}{\partial\theta_{i_{1}}}p_{\theta}(x)\cdots\frac{\partial}{\partial\theta_{i_{k}}}p_{\theta}(x)q_{\theta}(x)dx.
\end{align*}
\end{prop}
\begin{proof}
By Lemma \eqref{lem:biane}.
\end{proof}

~

\begin{prop}\label{pro:metric-converges}
The metric $g_{ij}(\theta,n)$ defined by equation \eqref{eq:metric}
converges.
\end{prop}
\begin{proof}
By definition $g_{ij}(\theta,n)=\frac{\partial^{2}}{\partial\theta_{i}\partial\theta_{j}}\psi(\theta,n)$.
\cite{ErMc} show that since $p_\theta(A)$ is a polynomial, the potential
$p_{\theta}(A)+\frac{1}{2}A^{2}$, has the expansion
\begin{align*}
\log\int\exp\left(-n\Tr\left(p_{\theta}(A)+\frac{1}{2}A^{2}\right)\right)dA=n^{2}\epsilon_{0}(\theta)+\epsilon_{1}(\theta)+\frac{1}{n^{2}}\epsilon_{2}(\theta)+\ldots
\end{align*}
with $\epsilon_i$ an analytic function of $\theta$ for $i=0,1,\ldots$.
Notice that
\begin{align*}
\frac{\partial^2}{\partial\theta_i\partial\theta_j}\frac{1}{n^2}\log\int\exp\left(-n\Tr\left(p_{\theta}(A)+\frac{1}{2}A^{2}\right)\right)dA=\frac{\partial^{2}}{\partial\theta_{i}\partial\theta_{j}}\psi(\theta,n),
\end{align*}
because
$\frac{\partial}{\partial\theta_i}\left(-n\Tr\left(\frac{1}{2}A^2\right)\right)=0$.

~

Therefore,
\begin{align*}
\frac{\partial^{2}}{\partial\theta_{i}\partial\theta_{j}}\psi(\theta,n)=\frac{\partial^{2}}{\partial\theta_{i}\partial\theta_{j}}\epsilon_{0}(\theta)+\frac{1}{n^{2}}\cdot\frac{\partial^{2}}{\partial\theta_{i}\partial\theta_{j}}\epsilon_{1}(\theta)+\frac{1}{n^{4}}\cdot\frac{\partial^{2}}{\partial\theta_{i}\partial\theta_{j}}\epsilon_{2}(\theta)+\ldots
\end{align*}

Thus, we have
\begin{align*}
\lim_{n\rightarrow\infty}g_{ij}(\theta,n)=\lim_{n\rightarrow\infty}\frac{\partial^{2}}{\partial\theta_{i}\partial\theta_{j}}\psi(\theta,n)=\frac{\partial^{2}}{\partial\theta_{i}\partial\theta_{j}}\epsilon_{o}(\theta).
\end{align*}
\end{proof}

~

\begin{prop}
For $\alpha\in[-1,1]$, the $\left(\alpha\right)$-connections converge.
\end{prop}
\begin{proof}
We have
\begin{align}\label{eq:connection-summand}
\begin{split}
&\int\textrm{Tr}\left(\frac{\partial^{2}}{\partial\theta_{k}\partial\theta_{i}}p_{\theta,n}\right)\textrm{Tr}\left(\frac{\partial}{\partial\theta_{j}}p_{\theta,n}\right)d\mu_{\theta,n}(A)=\\
&\int\textrm{Tr}\left(\frac{\partial^{2}}{\partial\theta_{k}\partial\theta_{i}}\psi(\theta,n)\right)\textrm{Tr}\left(F_j+\frac{\partial}{\partial\theta_{j}}\psi(\theta,n)\right)d\mu_{\theta,n}(A)=\\
&n\frac{\partial^{2}}{\partial\theta_{k}\partial\theta_{i}}\psi\left(\theta,n\right)\cdot\left(\int\Tr\left(F_{i}\left(A\right)\right)d\mu_{\theta,n}(A)-\int\Tr\left(F_{i}\left(A\right)\right)d\mu_{\theta,n}(A)\right)=0
\end{split}
\end{align}
where the second equality is due to equation
\eqref{eq:derivative-of-pressure}.

~

By equation \eqref{eq:connection-summand},
\begin{align*} \Gamma_{ijk}^{(\alpha),n}=-\frac{1-\alpha}{2}\cdot
n\int\textrm{Tr}\left(\frac{\partial}{\partial\theta_{i}}p_{\theta,n}\right)\textrm{Tr}\left(\frac{\partial}{\partial\theta_{j}}p_{\theta,n}\right)\textrm{Tr}\left(\frac{\partial}{\partial\theta_{k}}p_{\theta,n}\right)d\mu_{\theta,n}(A)\textrm{.}\end{align*}
Writing out
$\frac{\partial}{\partial\theta_{k}}g_{ij}\left(\theta,n\right)$ and
using equation \eqref{eq:connection-summand} again, we see that
\begin{align*}
\Gamma_{ijk}^{(\alpha),n}=-\frac{1-\alpha}{2}\cdot\frac{\partial}{\partial\theta_{k}}g_{ij}\left(\theta,n\right).
\end{align*}
The argument in Proposition \ref{pro:metric-converges} using \cite{ErMc} shows that
\begin{align*}
\frac{\partial}{\partial\theta_{k}}g_{ij}\left(\theta,n\right)=\frac{\partial^{3}}{\partial\theta_{i}\partial\theta_{j}\partial\theta_{k}}\psi\left(\theta,n\right)\rightarrow\frac{\partial^{3}}{\partial\theta_{i}\partial\theta_{j}\partial\theta_{k}}\epsilon_{0}(\theta).
\end{align*}
Therefore,
\begin{align*}
\lim_{n\rightarrow\infty}\Gamma_{ijk}^{(\alpha),n}=-\frac{1-\alpha}{2}\cdot\frac{\partial^{3}}{\partial\theta_{i}\partial\theta_{j}\partial\theta_{k}}\epsilon_{0}(\theta).
\end{align*}
\end{proof}

~

\begin{prop}
The dual coordinate system converges.
\end{prop}
\begin{proof}
According to Lemma \eqref{lem:biane},
\begin{align*}
\lim_{n\rightarrow\infty}\eta_{i}=\lim_{n\rightarrow\infty}\frac{1}{n}\int\Tr\left(F_{i}(A)\right)d\mu_{\theta,n}(A)=\int{F}_{i}(x)q_{\theta}(x)dx,
\end{align*}
\end{proof}

~

\begin{thm}
The Legendre transform of $\psi$ converges.
\end{thm}
\begin{proof}
According to the discussion in section \ref{sec:Legendre-Transform},
the Legendre trasnform of $\psi\left(\theta,n\right)$ is
\begin{align*}
\varphi(\theta,n)&=\frac{1}{n}\int\Tr\left(p_{\theta,n}\left(A\right)\right)d\mu_{\theta,n}(A)-\frac{1}{n}\int\Tr\left(p\left(A\right)\right)d\mu_{\theta,n}(A)\\
&=\frac{1}{n}\int\Tr\left(\sum_{i=1}^m\theta_iF_i(A)\right)d\mu_{\theta,n}(A)+\psi(\theta,n)
\end{align*}
By Lemma \eqref{lem:biane},
\begin{align*}
\frac{1}{n}\int\Tr\left(\sum_{i=1}^{m}\theta_{i}F_{i}(A)\right)d\mu_{\theta,n}(A)\rightarrow\int\left(\sum_{i=1}^{m}\theta_{i}F_{i}(x)\right)q_{\theta}(x)dx,
\end{align*}
and
\begin{align*}
&\psi(\theta,n)=\frac{1}{n^{2}}\log\int\exp\left(-nTr\left(p(A)+\sum_{i=1}^{m}\theta_{i}F_{i}(A)\right)\right)dA\longrightarrow\\
&\iint\log\left|x-y\right|q_{\theta}(x)q_{\theta}(y)dxdy-\int\left(p(x)+\sum_{i=1}^m\theta_i{F}_i(x)\right)q_{\theta}(x)dx.
\end{align*}

~

Therefore, $\lim_{n\rightarrow\infty}\varphi(\theta,n)=$
\begin{align*}
\int\left(\sum_{i=1}^m\theta_iF_i(x)\right)q_{\theta}(x)dx&+\iint\log\left\vert{x-y}\right\vert{q}_\theta(x)q_\theta(y)dxdy-\int\left(p(x)+\sum_{i=1}^m\theta_iF_i(x)\right)q_\theta(x)dx\\
&=\iint\log\left|x-y\right|q_{\theta}(x)q_{\theta}(y)dxdy-\int{p}(x)q_{\theta}(x)dx.
\end{align*}
\end{proof}

~

Let $q_\theta$ denote the limit of $d\mu_{\theta,n}$, and let $\chi(q_\theta)$ denote the Free Entropy of $q_\theta$. Since $\chi(q_\theta)=-\iint\log\vert{x-y}\vert{q}_\theta(x)q_\theta(y)dxdy$, the above calculation shows that
\begin{cor}\label{cor:entropy-converges}
\begin{align*}
H\left(d\mu_{\theta,n}\right)=-\frac{1}{n}\int\Tr\left(p_{\theta,n}(A)\right)d\mu_{\theta,n}(A)\rightarrow\chi\left(q_{\theta}\right)
\end{align*}
\end{cor}

Combining Corollary \ref{cor:entropy-converges} with Corollary
\ref{cor:convex} gives a new kind of convexity for free entropy of
the limit of a random matrix, which is "dual" to the convexity under
addition of the random variable:

~

\begin{cor}
Suppose $d\mu_{\theta,n}$ converges to $q_\theta$ for $\theta\in\Theta$ open. Then
$\chi\left(q_\theta\right)$ is convex in $\theta$.
\end{cor}

~

\begin{rem}
\cite{Hi} defines the free pressure for $R>0$ and $h\in{C}\left(\left[-R,R\right]\right)$ as
\begin{align*}
\pi_R\left(h\right)=\textrm{sup}\left\{\left.-\int{h}(x)d\mu(x)+\chi\left(\mu\right)\right|\mu\in\mathcal{M}\left(\left[-R,R\right]\right)\right\},
\end{align*}
where $\mathcal{M}([-R,R])$ is the set of Borel probability measures
supported on $[-R,R]$.  He calculates the free entropy as the
Legendre transform of free pressure with respect to a Banach space
duality,
\begin{align}\label{eq:hiai}
\chi\left(\mu\right)=\textrm{inf}\left\{\left.\int{h}(x)d\mu(x)+\pi_{R}\left(h\right)\right|h\in{C}\left(\left[-R,R\right]\right)\right\},
\end{align}

~

Now given a random matrix
$A\sim\exp\left(-n\textrm{Tr}\left(p(A)+\psi(n)\right)\right)$ with
$p\in\mathbb{R}\langle{x}\rangle$ convex, it converges to a measure $\mu$. Fix
$R$ so that $\textrm{supp}\left(\mu\right)\subset\left[-R,R\right]$,
and fix a collection $F_1,\ldots,F_k\in\mathbb{R}\langle{x}\rangle$. Consider the random matrix model
\begin{align*}
\exp\left(-n\textrm{Tr}\left(p(A)+\sum_{i=1}^m\theta_iF_i(A)+\psi(\theta,n)\right)\right).
\end{align*}
By definition of $\psi$ and $\pi_R$,
\begin{align*}
\psi(\theta,n)\rightarrow\pi_R\left(\sum_{i=1}^m\theta_iF_i(A)\right),
\end{align*}
and by Corollary \ref{cor:entropy-converges}
\begin{align*}
\varphi(0,n)\rightarrow \chi(\mu).
\end{align*}

We showed in section \ref{sec:Legendre-Transform} that
\begin{align*}
\varphi(0,n)=\inf_{\theta^{\prime}}\left\{\frac{1}{n}\int\textrm{Tr}\left(\sum_{i=1}^m\theta^{\prime}_iF_i(A)\right)d\mu_{0,n}(A)+\psi(\theta^{\prime},n)\right\},
\end{align*}
where the Legendre transform comes from a duality with respect to
$g_{ij}(\theta)$. In the limit $n\rightarrow\infty$ this is the
equation
\begin{align*}
\chi(\mu)=\textrm{inf}\left\{ \left.\int{h}(x)d\mu(x)+\pi_{R}\left(\mu\right)\right|h\in\textrm{span}\left\{F_i\right\}\right\},
\end{align*}
which is a restriction of \eqref{eq:hiai} to $h\in\textrm{span}\left\{F_i\right\}$.

~

Therefore, the restriction of Hiai's Banach space duality to any
finite linear span agrees with the corresponding Fisher information
metric duality.
\end{rem}

~

\begin{rem}
Given a noncommutative random variable $X\in(\mathcal{A},\tau)$ satisfying $\partial_X^*(1\otimes1)\in\poly{R}{x}$ and $\chi(X)<\infty$, we can uniquely define its information geometry (up to a constant) as follows.

~

Suppose $p\in\poly{R}{x}$ convex such that the random matrix model $\exp\left(-n\Tr\left(p(A)+\psi(n)\right)\right)$ converges to $X$ in the sense that $\frac{1}{n}\int\Tr(f(A))\exp\left(-n\Tr\left(p(A)+\psi(n)\right)\right)dA\rightarrow\tau(f(X))$ for all $f\in{C}(\mathbb{R})$. Fix perturbation functions $F_1,\ldots,F_m\in\poly{R}{x}$.

\begin{defn}
The \emph{information geometry of $X$ relative to $F_1,\ldots,F_m$} is the limit of the information geometry of $\exp\left(-n\Tr\left(p(A)+\sum_{i=1}^m\theta_iF_i+\psi(\theta,n)\right)\right)$ at $\theta=0$.
\end{defn}

To show uniqueness, suppose we have $q\in\poly{R}{x}$ convex such that the random matrix model\\
$\exp\left(-n\Tr\left(q(A)+\phi(n)\right)\right)$ also converges to $X$. By Lemma \ref{lem:biane} $p^\prime(x)=\int\frac{d\mu_X(y)}{x-y}=q^\prime(x)$, so $q(x)=p(x)+c$.
The constant $c$ may be absorbed into $\phi=\phi+c$. Since the tangent vectors, metric, and connections only depend on $\frac{\partial}{\partial\theta_i}\phi=\frac{\partial}{\partial\theta_i}\psi$ and $\frac{\partial^2}{\partial\theta_i\partial\theta_j}\phi=\frac{\partial^2}{\partial\theta_i\partial\theta_j}\psi$, the information geometries are the same at $\theta=0$.

To show existence, since $\chi(X)=-\iint\log\frac{1}{\vert{x-y}\vert}d\mu_X(x)d\mu_X(y)<\infty$, $p(x)=\int_{-\infty}^x\log\frac{1}{\vert{x-y}\vert}d\mu_X(x)d\mu_X(y)$ is well-defined and continuous. Since $p^{\prime\prime}(x)=2\int\frac{d\mu_X(y)}{\vert{x-y}\vert^2}>0$, $p$ is convex. By Lemma \eqref{lem:biane} the random matrix model $\exp\left(-n\Tr\left(p(A)+\psi(n)\right)\right)$ converges to $X$, and $p^\prime(x)\in\poly{R}{x}$. Thus, $p\in\poly{R}{x}$ is convex and its random matrix model converges to $X$.
\end{rem}

~

\subsection{Conjugate Variable and Free Fisher Information Measure}\label{sec:Conjugate-Variable}

~

One of the motivations for this paper was to understand
Voiculescu's conjugate variable. In this section we show that given
a random matrix model $\exp\left(-n\textrm{Tr}\left(p(A)+\psi(n)\right)\right)$,
converging to an operator on a Hilbert space $X$, we can construct a
random matrix model for $X$ with the tangent vector converging in
moments to the conjugate variable
$\partial_{X}^{*}\left(1\otimes1\right)$, and the Fisher information
metric converging to Voiculescu's Fisher information measure
$\Phi(X)$. Then we note that an analogous result holds for freely
independent $X_{1},\ldots,X_{k}$ using \cite[Prop 3.6]{Voi2}.

~

The calculation in Section \ref{sec:The-Tangent-Space} shows that
\begin{align*}
\partial_X^*(1\otimes 1)=-2\int\frac{d\mu_{X}(x)}{x-y}=p^\prime(y),
\end{align*}
which suggests the random matrix model
\begin{align*}
\big\{\exp\left(-n\Tr\left(p(A)+rp^\prime(A)+\psi(r,n)\right)\right)\big\}
\end{align*}
for the conjugate variable, with $\psi(r,n)=\frac{1}{n^{2}}\log\int\exp\left(-nTr\left(p(A)+rp^\prime(A)\right)\right)dA$.

The tangent vector to this model at $r=0$ is $p^{\prime}(A)+\left.\frac{\partial}{\partial r}\right|_{r=0}\psi(r,n)$, and the Fisher information metric at $r=0$ is
\begin{align*}
g_{11}(0)=\int\textrm{Tr}\left(p^{\prime}\left(A\right)+\left.\frac{\partial}{\partial{r}}\right|_{r=0}\psi(r,n)\right)\textrm{Tr}\left(p^{\prime}\left(A\right)+\left.\frac{\partial}{\partial{r}}\right|_{r=0}\psi(r,n)\right)d\mu_{0,n}(A).
\end{align*}

~

To evaluate these expressions we need a few formulae.

\begin{prop}\label{pro:integration-by-parts}
For $h\in\mathbb{R}\langle{x}\rangle$,
\begin{align*}
\int\Tr\left(h^{\prime}\left(A\right)\right)d\mu_{0,n}(A)=
n\int\Tr\left(h\left(A\right)\right)\Tr\left(p^{\prime}\left(A\right)\right)d\mu_{0,n}(A).
\end{align*}
\end{prop}
\begin{proof}
For a monomial, $h(A)=A^{k}$, we have
\begin{align*}
&\int\Tr\left(h^{\prime}\left(A\right)\right)d\mu_{0,n}(A)=\int\Tr(kA^{k-1})d\mu_{0,n}(A)=\\
&k\sum_{i_{1},\ldots,i_{k-1}}\int A_{i_{1}i_{2}}\cdots{A}_{i_{k-1} i_1}d\mu_{0,n}(A)=\\
&\sum_{j=1}^k\sum_{\substack{i_1,\ldots,i_k\\i_j=i_{j+1}}}\int A_{i_{1}i_{2}}\cdots{A}_{i_{j-1}i_j}\widehat{A_{i_ji_j+1}}A_{i_{j+1}i_{j+2}}\cdots A_{i_ki_1}d\mu_{0,n}(A)=\\
&\int\sum_{x=1}^{n}\sum_{i_{1},\ldots,i_{k}}\frac{\partial}{\partial{A}_{xx}}A_{i_{1}i_{2}}\cdots{A}_{i_ki_1}d\mu_{0,n}(A)=\sum_{x=1}^n\int\left(\frac{\partial}{\partial{A}_{xx}}\Tr\left(A^k\right)\right)d\mu_{0,n}=\\
&-\sum_{x=1}^n\left(-n\int\Tr\left(h\left(A\right)\right)\left(\frac{\partial}{\partial{A}_{xx}}\textrm{Tr}\left(p_{0,n}\left(A\right)\right)\right)\exp\left(-n\Tr\left(p_{0,n}(A)\right)\right)dA\right),
\end{align*}
with the last equality due to integration by parts.

~

Applying this calculation in reverse to $\frac{\partial}{\partial
A_{xx}}\Tr(p_{0,n})$, we see that $\frac{\partial}{\partial
A_{xx}}\Tr(p_{0,n})=\Tr(p^{\prime}_{0,n}(A))$.  So we arrive at
\begin{align*}
\int\Tr(h^\prime(A))d\mu_{0,n}(A)=n\int\Tr(h(A))\Tr(p_{0,n}^{\prime}(A))d\mu_{0,n}(A).
\end{align*}
By linearity of $\Tr$ we have the result for any $h\in\mathbb{R}\langle{x}\rangle$.
\end{proof}

~

Now we can evaluate using Proposition \ref{pro:integration-by-parts}
\begin{align*}
\left.\frac{\partial}{\partial{r}}\right|_{r=0}\psi(r,n)&=-\frac{1}{n}\int\Tr\left(p^{\prime}(A)\right)d\mu_{0,n}(A)=-\frac{1}{n^2}\int\Tr\left(1\right)\Tr\left(p^{\prime}(A)\right)d\mu_{0,n}(A)\\
&=-\frac{1}{n}\int Tr\left(0\right)d\mu_{0,n}(A)=0.
\end{align*}
So the tangent vector at $r=0$ is in fact $p^{\prime}(A)$. By Lemma
\eqref{lem:biane}, as $n\rightarrow\infty$
\begin{align*}
\frac{1}{n}\int\Tr\left(\left(p^{\prime}\left(A\right)\right)^{k}\right)d\mu_{0,n}(A)\rightarrow\int\left(p^{\prime}\left(x\right)\right)^{k}d\mu_{X}\left(x\right)=\int\left(\partial_{X}^{*}\left(1\otimes1\right)\left(x\right)\right)^{k}d\mu_{X}(x)
\end{align*}
so the tangent vector at $r=0$ indeed converges in moments to the
conjugate variable.

~

Now the metric at $r=0$ becomes
\begin{align*}
g_{11}\left(0,n\right)=\int\textrm{Tr}\left(p^{\prime}\left(A\right)\right)\textrm{Tr}\left(p^{\prime}\left(A\right)\right)d\mu_{0,n}(A).
\end{align*}

~

\begin{prop}
\label{pro: Fisher-Metric-Fisher-Measure}
\begin{align*}-\int\textrm{Tr}\left(p^{\prime}\left(A\right)\right)\textrm{Tr}\left(p^{\prime}\left(A\right)\right)d\mu_{0,n}(A)
 \rightarrow \Phi(X).\end{align*}
\end{prop}

\begin{proof}
We need a formula of \cite[Formula 2.18]{Jo}: if $\varphi\in
C^{1}(\mathbb{R})$ with $\varphi^{\prime}$ bounded below,

\begin{equation}\label{eq:johansson}
n\left(n-1\right)\int\frac{\varphi\left(t\right)-\varphi\left(s\right)}{t-s}u_{2,n}\left(s,t\right)dtds-n^{2}\int
p^{\prime}\left(t\right)\varphi\left(t\right)u_{1,n}\left(t\right)dt+n\int\varphi^{\prime}\left(t\right)u_{1,n}\left(t\right)dt
= 0\textrm{,}
\end{equation}

where we took $N,M=n$, $\beta=2$, and $h=0$ in the formula, and
\begin{align*}
u_{2,n}\left(\lambda_{1},\lambda_{2}\right)&=\int\prod_{a<b}\left(\lambda_{a}-\lambda_{b}\right)^{2}\prod_{c=1}^{n}\exp\left(-np\left(\lambda_{c}\right)\right)d\lambda_{3}\ldots{d}\lambda_{n},
\end{align*}
\begin{align*}
u_{1,n}\left(\lambda_{1}\right)&=\int\prod_{a<b}\left(\lambda_{a}-\lambda_{b}\right)^{2}\prod_{c=1}^{n}\exp\left(-np\left(\lambda_{c}\right)\right)d\lambda_{2}\ldots{d}\lambda_{n}.
\end{align*}

Notice that $d\mu_X(s)=\lim_{n\rightarrow\infty}u_{1,n}(s)$, where
$d\mu_X$ is the spectral measure of $X$. Now recall that
\begin{align*}
n\int\varphi^{\prime}\left(t\right)u_{1,n}\left(t\right)dt=\int\textrm{Tr}\left(\varphi^{\prime}(A)\right)d\mu_{0,n}(A)=n\int\textrm{Tr}\left(\varphi(A)\right)\textrm{Tr}\left(p^{\prime}(A)\right)d\mu_{0,n}(A)
\end{align*}
by Proposition \ref{pro:integration-by-parts}. Plugging this into
equation \eqref{eq:johansson}, and setting $\varphi=p^{\prime}$, we
get
\begin{align*}
&-\int\textrm{Tr}\left(p^{\prime}(A)\right)\textrm{Tr}\left(p^{\prime}(A)\right)d\mu_{0,n}(A)=n\int{p}^{\prime}\left(t\right)p^{\prime}\left(t\right)u_{1,n}\left(t\right)dt-\left(n-1\right)\int\frac{p^{\prime}\left(t\right)-p^{\prime}\left(s\right)}{t-s}u_{2,n}\left(s,t\right)dtds\\
&=\int{p}^{\prime}\left(t\right)p^{\prime}\left(t\right)u_{1,n}\left(t\right)dt+\left(n-1\right)\left(\int{p}^{\prime}\left(t\right)p^{\prime}\left(t\right)u_{1,n}\left(t\right)dt-\int\frac{p^{\prime}\left(t\right)-p^{\prime}\left(s\right)}{t-s}u_{2,n}\left(s,t\right)dtds\right).
\end{align*}

~

If we can show that
\begin{align}\label{eq:desired-johansson-formula}
\left(n-1\right)\left(\int{p}^{\prime}\left(t\right)p^{\prime}\left(t\right)u_{1,n}\left(t\right)dt-\int\frac{p^{\prime}\left(t\right)-p^{\prime}\left(s\right)}{t-s}u_{2,n}\left(s,t\right)dtds\right)\rightarrow0
\end{align}
then we would have
\begin{align*}
&-\int\textrm{Tr}\left(p^{\prime}(A)\right)\textrm{Tr}\left(p^{\prime}(A)\right)d\mu_{0,n}(A)=\int{p}^\prime\left(t\right)p^{\prime}\left(t\right)u_{1,n}\left(t\right)dt\rightarrow\\
&\int\left(p^{\prime}(t)\right)^2{d}\mu_X(t)=\int\left(\partial_X^{*}(1\otimes1)(t)\right)^2{d}\mu_X(t)=\Phi\left(X\right)
\end{align*}
and we would be done.

~

We use another fact from \cite[Prop. 2.6]{Jo} that for any
$\varphi\in C(\mathbb{R}^{2})$,
\begin{align*}
\lim_{n\rightarrow\infty}\int\varphi(s,t)u_{2,n}\left(s,t\right)dtds=\int\varphi(s,t)d\mu_{X}\left(s\right)d\mu_{X}\left(t\right).
\end{align*}
Thus, we have
\begin{align*}
&\lim_{n\rightarrow\infty}\int\frac{p^{\prime}\left(t\right)-p^{\prime}\left(s\right)}{t-s}u_{2,n}\left(s,t\right)dtds=-2\int{p}^{\prime}(s)\left(\int\frac{d\mu_{X}(t)}{t-s}\right)d\mu_{X}(s)=\\
&\int{p}^{\prime}(s)p^{\prime}(s)d\mu_{X}(s)=\lim_{n\rightarrow\infty}\int{p}^{\prime}\left(t\right)p^{\prime}\left(t\right)u_{1,n}\left(t\right)dt.
\end{align*}

We add and subtract this limit to equation
\eqref{eq:desired-johansson-formula} to get

\begin{align}\label{eq:convergence-1}
(n-1)\left(\int{p}^{\prime}\left(t\right)p^{\prime}\left(t\right)u_{1,n}\left(t\right)dt-\lim_{n\rightarrow\infty}\int{p}^{\prime}\left(t\right)p^{\prime}\left(t\right)u_{1,n}\left(t\right)dt\right)-
\end{align}
\begin{align}\label{eq:convergence-2}
(n-1)\left(\int\frac{p^{\prime}\left(t\right)-p^{\prime}\left(s\right)}{t-s}u_{2,n}\left(s,t\right)dtds-\lim_{n\rightarrow\infty}\int\frac{p^{\prime}\left(t\right)-p^{\prime}\left(s\right)}{t-s}u_{2,n}\left(s,t\right)dtds\right).
\end{align}

Now for \eqref{eq:convergence-1}, notice that
\begin{align*}
\int{p}^{\prime}\left(t\right)p^{\prime}\left(t\right)u_{1,n}\left(t\right)dt=\frac{1}{n}\int\textrm{Tr}\left(p^{\prime}(A)^{2}\right)d\mu_{0,n}(A)=\left.\frac{\partial}{\partial{s}}\right|_{s=0}\nu(s,n),
\end{align*}
with $\nu(s,n)=\frac{1}{n^{2}}\log\int\exp\left(-nTr\left(p(A)+s\cdot\left(p^{\prime}(A)\right)^{2}\right)\right)dA\textrm{.}$
Using the expansion of \cite{ErMc} we get
\begin{align*}\nu(s,n)=\epsilon_{o}(s)+\frac{1}{n^{2}}\epsilon_{1}(s)+\ldots\end{align*}
with $\epsilon_{i}(s)$ analytic. Therefore,
\begin{align*}\lim_{n\rightarrow\infty}\int
p^{\prime}\left(t\right)p^{\prime}\left(t\right)u_{1,n}\left(t\right)dt=\lim_{n\rightarrow\infty}\left.\frac{\partial}{\partial
s}\right|_{s=0}\psi(s,n)=\epsilon_{o}(s)\textrm{,}\end{align*}
and
\begin{align*}
\eqref{eq:convergence-1}=n\left(\left.\frac{\partial}{\partial{s}}\right|_{s=0}\nu(s,n)-\left.\frac{\partial}{\partial{s}}\right|_{s=0}\epsilon_{0}(s)\right)=n\left(\frac{1}{n^{2}}\left.\frac{\partial}{\partial{s}}\right|_{s=0}\epsilon_{1}(s)+\frac{1}{n^{4}}\left.\frac{\partial}{\partial{s}}\right|_{s=0}\epsilon_{2}(s)+\ldots\right)\longrightarrow0.
\end{align*}

~

For (\ref{eq:convergence-2}), we use \cite[formula 2.5]{Eyn}

\begin{align}
\begin{split}\label{eq:u_2}
u_{2,n}(s,t) = &\frac{n}{n-1}\left(\frac{1}{n}\sum_{l=0}^{n-1}P_{l}(s)^{2}\exp\left(-np(s)\right)\right)\left(\frac{1}{n}\sum_{l=0}^{n-1}P_{l}(t)^{2}\exp\left(-np(t)\right)\right)-\\
&\frac{n}{n-1}
\left(\alpha_{n}\cdot\frac{1}{n}\cdot\frac{P_{n}(s)P_{n-1}(t)-P_{n-1}(s)P_{n}(t)}{s-t}\exp\left(-\frac{n}{2}\left(p(t)+p(s)\right)\right)\right)^{2}.
\end{split}
\end{align}

where $P_{0},\ldots,P_{n}$ are the monic orthogonal polynomials with
respect to the measure $\exp\left(-np(t)\right)dt$, and $\alpha_{n}$
are constants converging to  a constant $\alpha$ which only depends
on the support of the limit distribution. \cite[after formula
2.5]{Eyn} gives the expansion
\begin{align*}
&P_{n}(s)=\frac{1}{\sqrt{f(s)}}\cos\left(n\xi(s)+g(s)\right)\cdot\exp\left(\frac{n}{2}p(s)\right),\\
&P_{n-1}(s)=\frac{1}{\sqrt{f(s)}}\cos\left(n\xi(s)+\varphi(s)+g(s)\right)\cdot\exp\left(\frac{n-1}{2}p(s)\right),
\end{align*}
where $f(s)$, $\xi(s)$, $\varphi(s)$, and $g(s)$ are functions of
$s$. We are only concerned with the order $n$ expansion, so we do
not need all these functions explicitly, but Eynard notes (in
formula (2.11)) that $f(s)=\sqrt{\left(s-a\right)\left(b-s\right)}$.
Recognizing the first summand in (\ref{eq:u_2}) as
$\frac{n}{n-1}\cdot u_{1,n}(s)\cdot u_{1,n}(t)$ and rewriting the
second summand using the expansion, we have
\begin{gather*}
u_{2,n}(s,t)=\frac{n}{n-1}\cdot{u}_{1,n}(s)\cdot{u}_{1,n}(t)-\frac{n}{n-1}\left(\alpha_{n}\cdot\frac{1}{n}\cdot\frac{1}{\sqrt{f(s)}\sqrt{f(t)}}\right)^{2}\times\\
\left(\frac{\cos\left(n\xi(s)+g(s)\right)\cos\left(n\xi(t)+\varphi(t)+g(t)\right)-\cos\left(n\xi(t)+g(t)\right)\cos\left(n\xi(s)+\varphi(s)+g(s)\right)}{s-t}\right)^{2}.
\end{gather*}

Now $-1\leq\cos\left(\cdot\right)\leq1$, so for $s,t$ satisfying
$\left|s-t\right|>n^{-1/8}$, $\left|s-a\right|>n^{-1/8}$ and
$\left|s-b\right|>n^{-1/8}$, the absolute value of the second
summand is bounded by
\begin{align*}
\frac{n}{n-1}\cdot\alpha_n^2\cdot n^{-2} \cdot n^{1/4}
\cdot n^{1/4} \longrightarrow 0.
\end{align*}

~

The first summand converges to $d\mu_{X}(s)d\mu_{X}(t)$, so we see
that the limit of $u_{2,n}(s,t)$ is $d\mu_{X}(s)d\mu_{X}(t)$ a.e.
$dsdt$. In equation \eqref{eq:convergence-2} we have subtracted
$\lim u_{2,n}(s,t)$, so we have
\begin{align*}
\eqref{eq:convergence-2}=&(n-1)\int\frac{p^{\prime}(t)-p^{\prime}(s)}{t-s}\left(\frac{n}{n-1}u_{1,n}(s)u_{1,n}(t)dsdt-\lim_{n\rightarrow\infty}u_{1,n}(s)u_{1,n}(t)dsdt\right)+\\
&(n-1)\int\frac{p^{\prime}(t)-p^{\prime}(s)}{t-s}\cdot\frac{n}{n-1}\left(\alpha_{n}\cdot\frac{1}{n}\cdot\frac{1}{\sqrt{f(s)}\sqrt{f(t)}}\right)^{2}\times\\
&\left(\frac{\cos\left(n\xi(s)+\chi(s)\right)\cos\left(n\xi(t)+\varphi(t)+\chi(t)\right)-\cos\left(n\xi(t)+\chi(t)\right)\cos\left(n\xi(s)+\varphi(s)+\chi(s)\right)}{s-t}\right)^{2}.
\end{align*}

~

In our calculation that $\eqref{eq:convergence-1}\rightarrow 0$, we
showed that $n\left(u_{1,n}-\lim u_{1,n}\right)\longrightarrow 0$,
so
\begin{align*}n\left(u_{1,n}(s)u_{1,n}(t)-\lim u_{1,n}(s)u_{1,n}(t)\right)\rightarrow 0.\end{align*}
Thus, the first summand in $\eqref{eq:convergence-2}\rightarrow 0$.

~

For $s,t$ satisfying $\left|s-t\right|>n^{-1/8}$,
$\left|s-a\right|>n^{-1/8}$ and $\left|s-b\right|>n^{-1/8}$,
\begin{gather*}
\left\vert\frac{n}{n-1}\left(\alpha_{n}\cdot\frac{1}{n}\cdot\frac{1}{\sqrt{f(s)}\sqrt{f(t)}}\right)^{2}\times\right.\\
\left.\left(\frac{\cos\left(n\xi(s)+\chi(s)\right)\cos\left(n\xi(t)+\varphi(t)+\chi(t)\right)-\cos\left(n\xi(t)+\chi(t)\right)\cos\left(n\xi(s)+\varphi(s)+\chi(s)\right)}{s-t}\right)^{2}\right\vert\leq\\
\frac{n}{n-1}\cdot\alpha_n^2\cdot n^{-1} \cdot n^{1/4}\cdot{n}^{1/4}\longrightarrow0.
\end{gather*}
Therefore, $\eqref{eq:convergence-2}\rightarrow0$. So
$\eqref{eq:desired-johansson-formula}\rightarrow0$, and we have proved the proposition.
\end{proof}

~

\begin{rem}
By \cite[Prop 3.6]{Voi2}, given $X_{1},\ldots,X_{k}$ freely
independent, the conjugate variable
$\partial_{X_i}^*(1\otimes1)$ computed in larger algebra
$L^2\left(\mathbb{R}\langle{X_1,\ldots,X_k}\rangle,\tau\right)$ satisfies
$\partial_{X_i}^*(1\otimes1)\in{L}^2(\mathbb{R}\langle{X_i}\rangle,\tau\vert_{\mathbb{R}\langle{X_i}\rangle})$.
Thus, $\partial_{X_i}^*(1\otimes1)=\frac{1}{2}\int\frac{d\mu_{X_i}(y)}{X_i-y}$. Thus, given $X_1,\ldots,X_k$ freely independent, with $X_i$ the limit of the random matrix model
$\exp\left(-nTr\left(p_i(A)\right)\right)$, we consider the following independent multi-matrix model (discussed in Section \ref{sec:Independent-Observations})
\begin{align*}
\exp\left(-n\Tr\left(\sum_{i=1}^{k}p_{i}(A_{i})+\sum_{j=1}^{k}t_{j}\cdot{p}_j^\prime(A_j)+\psi\left(t,n\right)\right)\right).
\end{align*}

Since the random matrices $A_{1},\ldots,A_{k}$ are independent, they
are asymptotically free, and since $A_{i}$ is distributed according
to $\exp\left(-nTr\left(p_{i}(A)\right)\right)$, the joint
distribution of $\left(A_{1},\ldots,A_{k}\right)$ converges to the
joint distribution of $\left(X_{1},\ldots,X_{k}\right)$. Also, the
tangent vector $p_{j}^{\prime}(A_{j})$ converge to
$p_{j}^{\prime}(X_{j})$ which is $\partial_{X_{j}}^{*}(1\otimes1)$.

According to the discussion in section
\ref{sec:Independent-Observations}, the off-diagonal terms in the metric vanish and we have
\begin{align*}
g_{jj}\left(0,n\right)=\int\Tr\left(p_j^\prime\left(A_j\right)\right)\Tr\left(p_j^\prime\left(A_{j}\right)\right)d\mu_{t_{j},n}(A_{j}).
\end{align*}
Applying Proposition \ref{pro: Fisher-Metric-Fisher-Measure}, we get
\begin{align*}
\int\Tr\left(p_j^\prime\left(A_j\right)\right)\Tr\left(p_j^\prime\left(A_j\right)\right)d\mu_{t_{j},n}(A_{j})\rightarrow\Phi\left(X_j\right).
\end{align*}
Therefore,
\begin{align*}
g_{ij}\left(0,n\right)\longrightarrow\left(
    \begin{array}{ccc}
    \Phi(X_1) &  & 0 \\
     & \ddots & \\
    0 &  & \Phi(X_k)
    \end{array}
\right).
\end{align*}
\end{rem}

~

~

\section{Cramer-Rao Theorem}\label{sec:Cramer-Rao-Theorem}

\subsection{Independent Observations}\label{sec:Independent-Observations}

~

In section (\ref{sec:Cramer-Rao-Theorem}) we will prove the
Cramer-Rao theorem, which requires us first to make sense of
independent observations and efficient estimators.

Classically, given a random variable $X$ whose distribution belongs
to a model\\ $S=\left\{\exp\left(p_{\theta}\left(x\right)\right)dx\left|\theta\in\Theta\subset\mathbb{R}^m\right.\right\}$, a common problem is to estimate the value of $\theta$ based on several independent observations of $X$.

Recall that a random variable is a real Borel function on some
probability space, so an observation of $X$ is simply a real number,
and $k$ observations of $X$ together form a vector in
$\mathbb{R}^k$. The requirement that the observations are
independent, means that $\left(x_1,\ldots,x_k\right)$ is an
observation of the random variable $X^{\otimes{k}}$ (shorthand
notation for $\big((X\otimes1\otimes\ldots\otimes1),\ldots,(1\otimes\ldots\otimes1\otimes{X})\big)$, which is
$k$ independent copies of $X$). Note that the distribution of
$X^{\otimes{k}}$ belongs to the model $S^{\oplus{k}}$, and the metric
on $S^{\oplus{k}}$ was calculated in Section \ref{sec:Several-Independent-Matrices} to be the direct sum of the metrics on each copy of $S$.

An \emph{estimator} is a collection of functions
$\xi_k:\mathbb{R}^k\rightarrow\mathbb{R}^m$, $k=1,2,\ldots$, used to estimated the value of $\theta\in\mathbb{R}^m$ based on $k$ independent observations.

An \emph{unbiased estimator} is an estimator such that if $X$ is
distributed according to $\exp\left(p_{\theta}\left(x\right)\right)dx$,
then the average of $\xi_{k}$ taken over all $k$ independent observations
of $X$ is equal to $\theta$. In concrete terms, independent observations
of $X$ means that $\theta\left(1\right)=\ldots=\theta\left(r\right)$,
so the requirement for an unbiased estimator becomes
\begin{align*}
\int\xi_k\left(x_1,\ldots,x_k\right)\exp\left(\sum_{r=1}^kp_{\theta}(x_r)\right)dx_1\ldots{d}x_k=\theta.
\end{align*}

Given an unbiased estimator, and $k$ independent observations
$\left(x_1,\ldots,x_k\right)$, the error of the estimate is
\begin{align*}
e_{\theta,k}\left(x_1,\ldots,x_k\right)=\xi_k\left(x_1,\ldots,x_k\right)-\theta.
\end{align*}

We may calculate the covariance matrix for the entries of the
error
\begin{align*}
\left(\textbf{Cov}\left(e_{\theta,k}\right)\right)_{ij}=\int\left(e_{\theta,k}\left(x_1,\ldots,x_k\right)\right)_{i}\left(e_{\theta,k}\left(x_{1},\ldots,x_{k}\right)\right)_{j}\exp\left(\sum_{r=1}^{k}p_{\theta}\left(x_{r}\right)\right)dx_{1}\ldots{d}x_k.
\end{align*}

The classical Cramer-Rao Theorem \cite{AmNa} gives a lower bound
on this covariance:

\begin{thm}
Given a model $S$ and an unbiased estimator $\xi_{k}$,
\begin{align*}
\textbf{Cov}\left(e_{\theta,k}\right)\geq{g}^{-1}\left(\theta\right)
\end{align*}
where $g$ is the Fisher Information metric on $S^{\oplus{k}}$,
and $\geq$ is in the sense of positive semi-definite matrices.
\end{thm}

~

In the random matrix case, we must first make sense of independent
observations. Recall that a random matrix $A$ is a matrix-valued
random variable, so

\begin{defn}
An observation of $A$ is simply a matrix $a\in{M}_n^{SA}(\mathbb{C})$, and $k$ observations of $A$ together form a vector $\left(a_1,\ldots,a_k\right)\in\left(M_n^{SA}(\mathbb{C})\right)^k$.
\end{defn}

\begin{defn}
$a_1,\ldots,a_k$ are $k$ \emph{independent observations} of $A$ if
$\left(a_1,\ldots,a_k\right)$ is an observation of $A^{\otimes{k}}$
($k$ independent copies of $A$).
\end{defn}

~

If the distribution of $A$ belongs to the model $S=\left\{\exp\left(-nTr\left(p_{n,\theta}(A)\right)\right)dA\right\}$, then the distribution of $A^{\otimes k}$ belongs to the
model $S^{\oplus{k}}$, which was discussed in section
(\ref{sec:Several-Independent-Matrices}). We denote
$P_{\theta,n}\left(A_{1},\ldots,A_{k}\right)$ for $p_{\theta,n}\left(A_{1}\right)+\ldots+p_{\theta,n}\left(A_{k}\right)$, and
$d\tilde{\mu}_{\theta,n}(A_{1},\ldots,A_{k})$ for $d\mu_{\theta,n}(A_{1})\cdots d\mu_{\theta,n}(A_{k})$.

~

The Cramer-Rao theorem also requires us to define unbiased
estimators. Recall from our previous discussion that an estimator is
a function on several independent copies of the random variable,
such that its expectation gives the estimated parameter value. In
addition, the proof of the classical Cramer-Rao rests on the fact
that an estimator may be viewed as a member of the tangent space of
the model. With these requirements in mind, we define

\begin{defn}
An \emph{estimator} is a collection of functions $(\xi_k)_i\in\mathbb{C}\langle{x_1,\ldots,x_k}\rangle$, where $k$ specifies the number of observations and $i$ specifies the parameter to be estimated. Given observations $A_{1},\ldots,A_{k}$, $\frac{1}{n}\textrm{Tr}\left(\xi_{k}\left(A_{1},\ldots,A_{k}\right)_{i}\right)$
is an estimate of $\theta_i$.
\end{defn}

\begin{defn}
An \emph{unbiased estimator} is an estimator such that
\begin{align*}
\frac{1}{n}\int\textrm{Tr}\left(\xi_{k}\left(A_{1},\ldots,A_{k}\right)_{i}\right)d\tilde{\mu}_{\theta,n}(A_{1},\ldots,A_{k})=\theta_{i}.
\end{align*}

\end{defn}
For an unbiased estimator and $k$ independent observations
$\left(A_{1},\ldots,A_{k}\right)$, we have the error of the $i^{th}$
estimate
\begin{align*}
\left(e_{k,\theta}\left(A_{1},\ldots,A_{k}\right)\right)_{i}=\left(\xi_{k}\left(A_{1},\ldots,A_{k}\right)-\theta\right)_{i},
\end{align*}
and we may compute the covariance matrix of the entries of the error using our inner-product
\begin{align*}
\textbf{Cov}\left(e_{k,\theta}\right)=-\int\textrm{Tr}\left(e_{k,\theta}\left(A_{1},\ldots,A_{k}\right)_{i}\right)\textrm{Tr}\left(e_{k,\theta}\left(A_{1},\ldots,A_{k}\right)_{j}\right)d\tilde{\mu}_{\theta,n}(A_{1},\ldots,A_{k}).
\end{align*}
We will bound this quantity in the Cramer-Rao theorem.

~

\subsection{Cramer-Rao Theorem}\label{sec:Cramer-Rao-Theorem}

~

In this section we prove the Cramer-Rao theorem for the random
matrix model. First we need a proposition:

\begin{prop}\label{pro:cramer-rao}
Given $f\in\polyxk$ a symmetric function, let $H:S^{\oplus{k}}\rightarrow\mathbb{R}$ be given by $\theta\mapsto-\frac{1}{n}\int\textrm{Tr}\left(f\left(A_{1},\ldots,A_{k}\right)\right)d\tilde{\mu}_{\theta,n}(A_{1},\ldots,A_{k})$. Then
\begin{align*}
\left\langle f+H(\theta),f+H(\theta)\right\rangle_{\theta}\geq\left\langle{dH,dH}\right\rangle_\theta,
\end{align*}
with equality if and only if $f\in T_{\theta}S^{\oplus{k}}$.
\end{prop}
\begin{proof}
First, consider the enlarged model
\begin{align*}
\tilde{S}=\left\{\exp\left(-n\Tr\left(\sum_{r=1}^{k}p\left(A_{r}\right)+\sum_{r=1}^{k}\sum_{i=1}^{m}\theta_{i}F_{i}\left(r\right)\left(A_{r}\right)+tf\left(A_{1},\ldots,A_{k}\right)+\psi\left(\theta,t,n\right)\right)\right)\right\}, \end{align*}
with $\psi(\theta,t,n)=\frac{1}{n^{2}}\log\int\exp\left(-n\textrm{Tr}\left(\sum_{r=1}^{k}p(A_{r})+\sum_{r=1}^{k}\sum_{i=1}^{m}\theta_{i}F_{i}\left(r\right)(A_{r})+tf(A_{1},\ldots,A_{k})\right)\right)dA_{1}\ldots{d}A_k$.
We have
\begin{align*}
\left.\left\langle{dH,dH}\right\rangle\right|_{\theta\in{S}^{\oplus{k+}}}\leq\left.\left\langle{dH,dH}\right\rangle\right|_{\left(\theta,0\right)\in\tilde{S}},
\end{align*}
with equality if and only if $\tilde{S}=S^{\oplus{k}}$. Let
$\tilde{\theta}=\left(\theta,t\right)$, and let
$P_{\tilde{\theta},n}\left(A_{1},\ldots,A_{k}\right)$ denote the
potential in $\tilde{S}$. For any $X\in{T}_{\left(\theta,0\right)}\tilde{S}$, we can differentiate $H$ with
respect to $X$ to get
\begin{align*}
X\cdot{H}(\theta)&=\sum_{i=1}^m\left\langle\frac{\partial{P}_{\tilde{\theta},n}}{\partial\tilde{\theta}_{i}},X\right\rangle_{\left(\theta,0\right)}\frac{\partial}{\partial\tilde{\theta}_{i}}H(\theta)=\sum_{i=1}^m\left\langle\frac{\partial{P}_{\tilde{\theta},n}}{\partial\tilde{\theta}_{i}},X\right\rangle_{\left(\theta,0\right)}\left\langle\frac{\partial{P}_{\tilde{\theta},n}}{\partial\tilde{\theta}_{i}},f\right\rangle_{\left(\theta,0\right)}\\
&=\sum_{i=1}^m\left\langle\frac{\partial{P}_{\tilde{\theta},n}}{\partial\tilde{\theta}_{i}},X\right\rangle_{\left(\theta,0\right)}\left\langle\frac{\partial{P}_{\tilde{\theta},n}}{\partial\tilde{\theta}_{i}},f+H(\theta)\right\rangle_{\left(\theta,0\right)}=\left\langle{f}+H(\theta),X\right\rangle_{\left(\theta,0\right)},
\end{align*}
where the next to last equality is due to the fact that
\begin{align*}
\left\langle\frac{\partial{P}_{\tilde{\theta},n}}{\partial\tilde{\theta}_{i}},H(\theta)\right\rangle=H(\theta)\frac{\partial}{\partial\tilde{\theta_{i}}}\int\exp\left(-nTr\left(P_{\tilde{\theta},n}\left(A_{1},\ldots,A_{k}\right)\right)\right)dA=0.
\end{align*}

~

The equation $X\cdot{H}(\theta)=\left\langle{f}+H(\theta),X\right\rangle_{\tilde{\theta}}$ is the definition of $\textrm{grad}\left(H\right)$ in $T_{\left(\theta,0\right)}\tilde{S}$, and combined with the fact that $f+H(\theta)\in T_{\left(\theta,0\right)}\tilde{S}$ it shows that
\begin{align*}
f+H(\theta)=\textrm{grad}\left(H\right)\textrm{ in
}T_{\left(\theta,0\right)}\tilde{S}.
\end{align*}

~

We have
\begin{align*}
&\left.\left\langle{dH,dH}\right\rangle\right|_{\theta\in{S}^{\oplus{k}}}\leq\left.\left\langle{dH,dH}\right\rangle\right|_{\left(\theta,0\right)\in\tilde{S}}=\left.\left\langle\textrm{grad}\left(H\right),\textrm{grad}\left(H\right)\right\rangle\right|_{\left(\theta,0\right)\in\tilde{S}}=\\
&\left.\left\langle{f}+H(\theta),f+H(\theta)\right\rangle\right|_{\left(\theta,0\right)\in\tilde{S}}=\left.\left\langle{f}+H(\theta),f+H(\theta)\right\rangle\right|_{\theta\in{S}}, \end{align*}
with equality if and only if $\tilde{S}=S^{\oplus{k}}$.

~

\end{proof}

~

Now we prove the Cramer-Rao theorem:

\begin{thm}\label{thm:cramer-rao}
Given an exponential family $S$, and an unbiased estimator $\xi_{k}$
on $S$, it satisfies \begin{align*} \left\langle
\left(\xi_{k}-\theta\right)_{i},\left(\xi_{k}-\theta\right)_{j}\right\rangle
_{\theta}\geq G^{-1}(\theta)\end{align*} in the sense of positive
semi-definite matrices, where $\left\langle \cdot,\cdot\right\rangle
_{\theta}$ is the inner-product on $S^{\oplus k}$, and $G$ is the
Fisher Information metric on $S^{\oplus k}$.
\end{thm}
\begin{proof}
To show that $\left(\left\langle\left(\xi_{k}-\theta\right)_{i},\left(\xi_{k}-\theta\right)_{j}\right\rangle_{\theta}\right)-g^{-1}$ is positive semi-definite, fix an arbitrary vector $c\in\mathbb{R}^m$, and let
\begin{align*}
f\left(A_{1},\ldots,A_{k}\right)=\sum_{i=1}^{m}c_{i}\left(\xi_{k}\left(A_{1},\ldots,A_{k}\right)\right)_{i}.
\end{align*}
Write $S^{\oplus{k}}=\left\{\exp\left(-nTr\left(\sum_{r=1}^{k}p_{\theta,n}\left(A_{r}\right)\right)\right)\right\}$, $P_{\theta,n}\left(A_{1},\ldots,A_{k}\right)$ for $\sum_{r=1}^{k}p_{\theta,n}\left(A_{r}\right)$, $d\tilde{\mu}_{\theta,n}(A_{1},\ldots,A_{k})$ for $d\mu_{\theta,n}(A_{1})\cdots{d}\mu_{\theta,n}(A_{k})$, and denote the expectation by
$\mathbb{E}_{\theta}\left(f\right)=\frac{1}{n}\int\textrm{Tr}\left(f\left(A_{1},\ldots,A_{k}\right)\right)d\tilde{\mu}_{\theta,n}(A_{1},\ldots,A_{k})$.
This model was discussed in Section \ref{sec:Several-Independent-Matrices}; denote its metric by $G_{ij}\left(\theta,n\right)$. Define $H(\theta)=-\mathbb{E}_{\theta}\left(f\right)$ as in
Proposition \ref{pro:cramer-rao}.

~

Because $\xi_k$ is unbiased, $H(\theta)=-\sum_{i=1}^kc_i\theta_i$, so
$\frac{\partial}{\partial\theta_i}H(\theta)=-c_i$, and Proposition \ref{pro:cramer-rao} shows that
\begin{align*}
\left\langle{dH,dH}\right\rangle_\theta\leq\left\langle{f+H}(\theta),f+H(\theta)\right\rangle_\theta.
\end{align*}

~

The left-hand side is
\begin{align*}
&\left\langle{dH,dH}\right\rangle_\theta=\left\langle\sum_{i=1}^{m}\frac{\partial}{\partial\theta_{i}}H(\theta)\cdot\frac{\partial{P}_{\theta,n}}{\partial\theta_{i}},\sum_{j=1}^{m}\frac{\partial}{\partial\theta_{j}}H(\theta)\cdot\frac{\partial{P}_{\theta,n}}{\partial\theta_{j}}\right\rangle_\theta=\\
&\sum_{i,j=1}^{m}\frac{\partial}{\partial\theta_{i}}\left(\sum_{r=1}^{m}c_{r}\mathbb{E}_{\theta}\left(\left(\xi_{k}\right)_{r}\right)\right)\frac{\partial}{\partial\theta_{j}}\left(\sum_{s=1}^{m}c_{s}\mathbb{E}_{\theta}\left(\left(\xi_{k}\left(A_{1},\ldots,A_{k}\right)\right)_{s}\right)\right)\left(G^{-1}\left(\theta,n\right)\right)_{ij}=\\
&\sum_{i,j=1}^{m}\frac{\partial}{\partial\theta_{i}}\left(\sum_{r=1}^{m}c_{r}\theta_{r}\right)\frac{\partial}{\partial\theta_{j}}\left(\sum_{s=1}^{m}c_{s}\theta_{s}\right)\left(G^{-1}\left(\theta,n\right)\right)_{ij}=\\
&\sum_{i,j=1}^{m}c_{i}c_{j}\left(G^{-1}\left(\theta,n\right)\right)_{ij}=c^{t}\left(G^{-1}\left(\theta,n\right)\right)_{ij}c,
\end{align*}
and the right-hand side is
\begin{align*}
\left\langle{f+H(\theta),f+H(\theta)}\right\rangle_\theta=\sum_{i,j=1}^{m}c_{i}c_{j}\left\langle\left(\xi_{k}-\theta\right)_{i},\left(\xi_{k}-\theta\right)_{j}\right\rangle_{\theta}.
\end{align*}
Therefore, we have shown that
\begin{align*}
c^{t}\left\langle\left(\xi_{k}-\theta\right)_{i},\left(\xi_{k}-\theta\right)_{j}\right\rangle_{\theta}c\geq{c}^{t}\left(G^{-1}\left(\theta,n\right)\right)_{ij}c,
\end{align*}
with equality if and only if $\left(\xi_k-\theta\right)_i\in{T}_\theta{S}^{\oplus{k}}$ for $i=1,\ldots,m$.

~

\end{proof}

~

Next we prove the converse to the Cramer-Rao theorem, which requires
a definition:

\begin{defn}
An \emph{efficient estimator} is an unbiased estimator $\xi_{k}$
that attains the bound in the Cramer-Rao theorem.
\end{defn}

~

\begin{thm}
A random matrix model $S=\left\{\exp\left(-nTr\left(Q_{\theta,n}(A)\right)\right)\right\}$, with $Q_{\theta}\in\poly{R}{x}$, has efficient estimators
$\hat{\theta}_{1},\ldots,\hat{\theta}_{k}$ if and only if it is an
exponential family, i.e.
\begin{align*}
Q_{\theta}\left(A\right)=p\left(A\right)+\sum\theta_{i}F_{i}\left(A\right)+\psi\left(\theta,n\right).
\end{align*}
\end{thm}
\begin{proof}
$\left(\Leftarrow\right)$ Suppose $S$ is an exponential family.
Recall that we defined the dual coordinate system
$\eta_{i}=\frac{\partial}{\partial\theta_{i}}\psi\left(\theta,n\right)$.
Thus, we may write the tangent vectors as $F_{i}-\eta_{i}$, and the
Fisher information metric as
\begin{align*}
g_{ij}\left(\theta,n\right) = \left\langle
\left(F_{i}-\eta_{i}\right),\left(F_{j}-\eta_{j}\right)\right\rangle
_{\theta}\textrm{.}\end{align*} This equation shows that \begin{align*}
\hat{\eta}_{i}\left(A_{1},\ldots,A_{k}\right) =
\frac{1}{k}\sum_{j=1}^{k}F_{i}\left(A_{j}\right)\end{align*}
 is an efficient estimator for $\eta_{i}$. Thus, we have found a
coordinate system $\eta_{1},\ldots,\eta_{k}$ for the exponential
family which has efficient estimators $\hat{\eta}_{1},\ldots,\hat{\eta}_{k}$.

$\left(\Rightarrow\right)$ Given $S$ with efficient estimators
$\hat{\theta}_{1}\left(A_{1},\ldots,A_{k}\right),\ldots,\hat{\theta}_{m}\left(A_{1},\ldots,A_{k}\right)$,
fix $k=1$. Since the estimators are efficient, they attain equality
in the Cramer-Rao theorem:
\begin{align*}
\left\langle\left(\hat{\theta}\left(A\right)-\theta\right)_{i},\left(\hat{\theta}\left(A\right)-\theta\right)_{j}\right\rangle_{\theta}=\left(G^{-1}\left(\theta,n\right)\right)_{ij},
\end{align*}
and this implies that $\underbrace{\hat{\theta}_{i}\left(A\right)-\theta_{i}}_{F_{i}\left(A\right)}\in{T}_{\theta}S$ by the equality in Proposition \ref{pro:cramer-rao}. So we have $m$ linearly independent vector fields on $S$. To see that they are parallel with respect to the
$(1)$-connection, fix $\theta^{\prime}$ and consider the
model
\begin{align*} S^{\prime} = \left\{
\exp\left(-n\textrm{Tr}\left(Q_{\theta^{\prime}}\left(A\right)+\sum_{i=1}^m\zeta_{i}F_{i}\left(A\right)+\phi\left(\zeta,n\right)\right)\right)\right\}
\textrm{,}\end{align*} with
$\phi\left(\zeta,n\right)=\frac{1}{n^{2}}\log\int\exp\left(-n\textrm{Tr}\left(Q_{\theta^{\prime}}\left(A\right)+\sum_{i=1}^m\zeta_{i}F_{i}\left(A\right)\right)\right)dA$.
Denote \begin{align*} p_{\zeta,n}(A) = Q_{\theta^{\prime}}(A)+\sum_{i=1}^m
\zeta_{i}F_{i}\left(A\right)+\phi\left(\zeta,n\right)\textrm{.}\end{align*}

~

This is a submanifold of the original model, and since
$F_{1},\ldots,F_{m}$ are linearly independent
$T_{\zeta=0}S^{\prime}=T_{\theta=\theta^{\prime}}S\textrm{.}$

Now
\begin{align*}
\left.\frac{\partial\phi}{\partial\zeta_{i}}\right|_{\zeta=0}&=-\frac{1}{n}\int\textrm{Tr}\left(F_{i}\left(A\right)\right)\exp\left(-n\textrm{Tr}\left(p_{\theta^{\prime},n}\left(A\right)\right)\right)dA\\
&=-\frac{1}{n}\int\textrm{Tr}\left(\hat{\theta}_{i}(A)-\theta_{i}\right)\exp\left(-n\textrm{Tr}\left(p_{\theta^{\prime},n}(A)\right)\right)dA=0,
\end{align*}
because $\hat{\theta}$ is an unbiased estimator. Thus,

\begin{align*}
\left.\frac{\partial{p}_{\zeta,n}}{\partial\zeta_{i}}\right|_{\zeta=0}&=-\frac{1}{n}\int\textrm{Tr}\left(F_{i}+\left.\frac{\partial\phi}{\partial\zeta_{i}}\right|_{\zeta=0}\right)\exp\left(-n\textrm{Tr}\left(p_{\theta^{\prime},n}\left(A\right)\right)\right)dA\\
&=-\frac{1}{n}\int\textrm{Tr}\left(F_{i}\left(A\right)\right)\exp\left(-n\textrm{Tr}\left(p_{\theta^{\prime},n}\left(A\right)\right)\right)dA=0
\end{align*}
again because $\hat{\theta}$ is an unbiased estimator. Now we calculate
the $(1)$-connection for the $\zeta$ coordinate system at $\zeta=0$:
\begin{align*}
\Gamma_{ijk}^{(-1),n}\left(0\right)&=\int\textrm{Tr}\left(\frac{\partial^{2}}{\partial\zeta_{k}\partial\zeta_{i}}p_{\zeta,n}\right)\textrm{Tr}\left(\frac{\partial}{\partial\zeta_{j}}p_{\zeta,n}\right)\exp\left(-n\textrm{Tr}\left(p_{\theta^{\prime},n}\left(A\right)\right)\right)dA\\
&=n\frac{\partial^{2}}{\partial\zeta_{k}\partial\zeta_{i}}\phi\cdot\int\textrm{Tr}\left(\frac{\partial}{\partial\zeta_{j}}p_{\zeta,n}\right)\exp\left(-n\textrm{Tr}\left(p_{\theta^{\prime},n}\left(A\right)\right)\right)dA=0.
\end{align*}
Therefore, $\zeta$ is a $(1)$-flat coordinate system for $S$ at
$\theta^{\prime}$ for $\theta^\prime\in{S}$, so we may write
\begin{align*}
Q_{\theta,n}\left(A\right)=Q_{\theta^{\prime},n}\left(A\right)+\sum\zeta_{i}\left(\hat{\theta}_{i}\left(A\right)-\theta_{i}\right)+\phi\left(\zeta,n\right).
\end{align*}
\end{proof}

~

\subsection{Relation to the Free Cramer-Rao}\label{sec:Relation-to-Free-Cramer-Rao}

~

In this section we compare the Cramer-Rao Theorem \ref{thm:cramer-rao} as $n\rightarrow\infty$ to Voiculescu's Free Cramer-Rao Theorems.
We recall the free Cramer-Rao theorem in the one-variable case
\cite{Voi2}:

\begin{thm}
Given $v\geq0$, $v\in L^{1}(\mathbb{R})\cap L^{3}(\mathbb{R})$,
$\int v(x)dx=1$, and $\int x^{2}v(x)dx<\infty$, let $x_{0}=\int xv(x)dx$,
then
\begin{align*}
\underbrace{\left(\int{v}(x)^{3}dx\right)}_{\Phi\left(v(x)dx\right)}\left(\int\left(x-x_{0}\right)^{2}v(x)dx\right)\geq\frac{3}{4\pi^{2}}.
\end{align*}

\end{thm}
In this theorem $f(x)=x$ plays the role of an estimator. To put this
in the conext of the random matrix model, consider a random matrix
$A\sim\exp\left(-nTr\left(p(A)\right)\right)dA$, with limit
distribution $v(x)$. In section (\ref{sec:Conjugate-Variable}) we
showed that Voiculescu's $\Phi$ corresponds to embedding $A$ in the
model\begin{align*}
\exp\left(-n\textrm{Tr}\left(p(A)+tp^{\prime}(A)+\psi(t,n)\right)\right)\textrm{,}\end{align*}
and considering the Fisher information metric at $t=0$. The
estimator is $f\left(A\right)=A$. To apply our Cramer-Rao theorem,
we must assume that $f$ is an unbiased estimator at $t=0$: \begin{align*}
\frac{1}{n}\int\textrm{Tr}\left(f\left(A\right)\right)\exp\left(-n\textrm{Tr}\left(p(A)\right)\right)dA=0\textrm{.}\end{align*}
 Our Cramer-Rao theorem says \begin{align*}
\left\langle f\left(A\right),f\left(A\right)\right\rangle _{t=0}\geq
G^{-1}(0,n)\textrm{,}\end{align*} where $G(t,n)$ is the metric on our model.
In section (\ref{sec:Conjugate-Variable}) we showed that
$G(0,n)\rightarrow\Phi(v(x)dx)$.

Thus, we have
\begin{align*}
\lim_{n\rightarrow\infty}\left(-\int\textrm{Tr}\left(A\right)\textrm{Tr}\left(A\right)\exp\left(-n\textrm{Tr}\left(p(A)\right)\right)dA\right)
 \geq \Phi^{-1}(v(x)dx)\textrm{.}\end{align*}

Now, since $\frac{1}{n}\textrm{Tr}\left(\left(A-\textrm{Tr}\left(A\right)\right)^{2}\right)=\frac{1}{n}\textrm{Tr}\left(A^{2}\right)-\frac{1}{n}\textrm{Tr}\left(A\right)\textrm{Tr}\left(A\right)+\textrm{Tr}\left(A\right)\textrm{Tr}\left(A\right)$, and $\frac{1}{n}\textrm{Tr}\left(\left(A-\textrm{Tr}\left(A\right)\right)^{2}\right)\geq0$, and $\frac{1}{n}\textrm{Tr}\left(A\right)\textrm{Tr}\left(A\right)\geq0$, we have
\begin{align*}
-\textrm{Tr}\left(A\right)\textrm{Tr}\left(A\right)\leq\frac{1}{n}\textrm{Tr}\left(A^{2}\right)\textrm{.}
\end{align*}
Thus, our Cramer-Rao implies that
\begin{align*}
\lim_{n\rightarrow\infty}\frac{1}{n}\int\textrm{Tr}\left(A^{2}\right)\exp\left(-n\textrm{Tr}\left(p(A)\right)\right)dA\geq\Phi^{-1}(v(x)dx)
\end{align*}
which is precisely
\begin{align*}
\int{x}^{2}v(x)dx\geq\int\left(v(x)\right)^{3}dx.
\end{align*}

~

~

The several variable case result \cite{Voi2} states that

\begin{thm}
Given $X_{1},\ldots,X_{k}\in\left(\mathcal{A},\tau\right)$ with $(\mathcal{A},\tau)$ a tracial unital von Neumann algebra,
\begin{align*}
\tau\left(X_{1}^{2}+\ldots+X_{k}^{2}\right)\geq k^{2}\cdot\Phi^{-1}\left(X_{1},\ldots,X_{k}\right)\textrm{.}\end{align*}

\end{thm}
To compare it to our theorem, we must assume that
$\left(X_{1},\ldots,X_{k}\right)$ are freely independent, and that
we have random matrices $A_{1},\ldots,A_{k}$, with $A_{r}$
converging to $X_{r}$, and $A_{r}$ distributed according to
$Z_n^-1\exp\left(-n\textrm{Tr}\left(p\left(r\right)\left(A\right)\right)\right)$.
Following section (\ref{sec:Several-Independent-Matrices}),
$\exp\left(-n\textrm{Tr}\left(\sum_{r=1}^{k}p\left(r\right)\left(A_{r}\right)\right)\right)$
converges to $\left(X_{1},\ldots,X_{k}\right)$. We embed this in the
model \begin{align*} S = \left\{
\exp\left(-n\textrm{Tr}\left(\sum_{r=1}^{k}p\left(r\right)\left(A_{r}\right)+t\sum_{r=1}^{k}p^{\prime}\left(r\right)\left(A_{r}\right)+\psi\left(t,n\right)\right)\right)\right\}
\textrm{,}\end{align*} and the discussion in section
(\ref{sec:Conjugate-Variable}) shows that
\begin{align*}
\Phi\left(X_{1},\ldots,X_{k}\right) =
\sum_{r=1}^{k}\lim_{n\rightarrow\infty}G_{rr}\left(0,n\right)\textrm{.}\end{align*}
The estimator in Voiculescu's theorem is
$f\left(A_{1},\ldots,A_{k}\right)=\frac{1}{k}\sum_{r=1}^{k}A_{r}$,
and we must assume it is unbiased at $t=0$:
\begin{align*}
\frac{1}{n}\int\textrm{Tr}\left(\frac{1}{k}\sum_{s=1}^{k}A_{s}\right)\exp\left(-n\textrm{Tr}\left(\sum_{r=1}^{k}p\left(r\right)\left(A_{r}\right)\right)\right)dA_{1}\ldots
dA_{k} = 0\textrm{.}\end{align*}

~

Now our Cramer-Rao theorem says
\begin{align*}
\int\textrm{Tr}\left(\frac{1}{k}\sum_{l=1}^{k}A_{l}\right)\textrm{Tr}\left(\frac{1}{k}\sum_{s=1}^{k}A_{s}\right)\exp\left(-n\textrm{Tr}\left(\sum_{r=1}^{k}p\left(r\right)\left(A_{r}\right)\right)\right)dA_{1}\ldots
dA_{k} \geq G^{-1}(0,n)\textrm{.}\end{align*}

~

~

Since $A_{l}$ is independent of $A_{s}$ for $l\neq s$, this is
equivalent to
\begin{align*}
\frac{1}{k^{2}}\sum_{s=1}^{k}\int\textrm{Tr}\left(A_{s}\right)\textrm{T}r\left(A_{s}\right)\exp\left(-n\textrm{Tr}\left(\sum_{r=1}^{k}p\left(r\right)\left(A_{r}\right)\right)\right)dA_{1}\ldots
dA_{k} \geq G^{-1}(0,n)\textrm{.}\end{align*}

This is equivalent to
\begin{align*}
\lim_{n\rightarrow\infty}\sum_{s=1}^{k}\int\Tr\left(A_{s}\right)\Tr\left(A_{s}\right)\exp\left(-n\Tr\left(\sum_{r=1}^{k}p\left(r\right)\left(A_{r}\right)\right)\right)dA_{1}\ldots{d}A_k\geq{k}^2\cdot\Phi^{-1}\left(X_{1},\ldots,X_{k}\right).
\end{align*}
Again, since
$-\textrm{Tr}\left(A_{s}\right)\textrm{Tr}\left(A_{s}\right)\leq\frac{1}{n}\Tr\left(A_{s}^{2}\right)$
for $s=1,\ldots,k$, we obtain
\begin{align*}
\lim_{n\rightarrow\infty}\sum_{s=1}^k\frac{1}{n}\int\Tr(A_s^2)\exp\left(-n\Tr\left(\sum_{r=1}^{k}p\left(r\right)\left(A_{r}\right)\right)\right)dA_{1}\ldots{d}A_k\geq{k}^2\cdot\Phi^{-1}\left(X_{1},\ldots,X_{k}\right)
\end{align*}
which is
\begin{align*}
\tau\left(X_{1}^{2}+,\ldots,X_{k}^{2}\right)\geq{k}^2\cdot\Phi^{-1}\left(X_{1},\ldots,X_{k}\right).
\end{align*}

~

\subsection{Relation to Second-Order Freeness}\label{sec:Second-Order-Freeness}

~

In this section we show that the quantities which motivated
Speicher's theory of second-order freeness are naturally related to
the geometry of the random matrix model.

Given a random matrix $A\sim\exp\left(-n\Tr\left(p(A)+\psi(n)\right)\right)$, and
a collection of functions $F_{1},\ldots,F_{m}\in\poly{R}{x}$,
\cite{MSp} provide several definitions for the fluctuations of
$F_{1}(A),\ldots,F_{m}(A)$, and we use one which is convenient for
us:
\begin{defn}
The fluctuations of $F_{1}\left(A\right),\ldots,F_{m}\left(A\right)$
are (if they exist)
\begin{align*}
\beta_{ij}=\lim_{n\rightarrow\infty}\int\Tr\left(F_{i}(A)-\alpha_{i}(n)\right)\Tr\left(F_{j}(A)-\alpha_{j}(n)\right)\exp\left(-n\textrm{Tr}\left(p(A)+\psi(n)\right)\right)dA,
\end{align*}
with
\begin{align*}
\alpha_{i}(n)=\frac{1}{n}\int\Tr\left(F_{i}(A)\right)\exp\left(-n\textrm{Tr}\left(p(A)+\psi(n)\right)\right)dA.
\end{align*}
\end{defn}

\cite{MSp} have developed a general framework of conditions on an
algebra so that the fluctuations of any of its elements may be
calculated; they call it Second-Order Freeness.

~

Consider the random matrix model
\begin{align*}
\exp\left(-n\textrm{Tr}\left(p(A)+\sum_{i=1}^m\theta_{i}F_{i}(A)+\psi(\theta,n)\right)\right).
\end{align*}
For small enough $\theta$, this model converges, and in particular
we have at $\theta=0$:
\begin{align*}
g_{ij}\left(0,n\right)=\int\textrm{Tr}\left(F_{i}+\left.\frac{\partial}{\partial\theta_{i}}\right|_{\theta=0}\psi\right)\textrm{Tr}\left(F_{j}+\left.\frac{\partial}{\partial\theta_{j}}\right|_{\theta=0}\psi\right)\exp\left(-n\textrm{Tr}\left(p(A)+\psi(0,n)\right)\right)dA.
\end{align*}
Since
\begin{align*}
\left.\frac{\partial}{\partial\theta_{i}}\right|_{\theta=0}\psi(\theta,n)=-\frac{1}{n}\int\textrm{Tr}\left(F_{i}\right)\exp\left(-n\textrm{Tr}\left(p(A)+\psi(0,n)\right)\right)dA=-\alpha_{i}(n),
\end{align*}
we have
\begin{align*}
g_{ij}\left(0,n\right)=\int\textrm{Tr}\left(F_{i}-\alpha_{i}(n)\right)\textrm{Tr}\left(F_{j}-\alpha_{j}(n)\right)\exp\left(-n\textrm{Tr}\left(p(A)+\psi(0,n)\right)\right)dA.
\end{align*}

~

In Section \ref{sec:The-n-goes-to-infinity-case} we showed that
this metric converges as $n\rightarrow\infty$, so we
have
\begin{align*}
g_{ij}(0)=\lim_{n\rightarrow\infty}\int\textrm{Tr}\left(F_{i}-\alpha_{i}(n)\right)\textrm{Tr}\left(F_{j}-\alpha_{j}(n)\right)\exp\left(-nTr\left(p(A)+\psi(0,n)\right)\right)dA.
\end{align*}
These are precisely the fluctuations of $F_{1}(A),\ldots,F_{m}(A)$.

~

\begin{rem}
As a result, we have shown that the fluctuations of a random matrix
may be considered as tangent vectors of the random matrix model
obtained by perturbing the potential with the fluctuation functions;
the inner-product of the fluctuations is the metric on this model; and as a result, the fluctuations of a random matrix give rise to a positive-definite form.
\end{rem}
~
\bibliography{Paper1}
\bibliographystyle{alpha}
~

\address{Department of Mathematics, University of California, Los Angeles
90095, USA}

\email{\emph{Email address:} dshiber@math.ucla.edu}
\end{document}